\documentclass[11pt]{article}
\usepackage{graphicx}
\usepackage{mathrsfs,amsmath,amssymb,amsfonts}
\numberwithin{equation}{section}
\usepackage{amsfonts}
\usepackage{amssymb, amsmath, cite}
\usepackage{color}
\usepackage{epsf}
\usepackage{epstopdf}
\usepackage{multirow}
\usepackage{bm}
\usepackage{booktabs}
\usepackage{float}

\newcommand{\ignore}[1]{}{}
\newcommand{\be}{\begin{equation}}               
\newcommand{\ee}{\end{equation}}                 
\newcommand{\bi}{\begin{itemize}}
\newcommand{\ei}{\end{itemize}}

\newcommand{\beqn}{\begin{eqnarray}}             
\newcommand{\eeqn}{\end{eqnarray}}               
\newcommand{\beq}{\begin{eqnarray*}}             
\newcommand{\eeq}{\end{eqnarray*}}               

\newcommand{\ssb}{\scriptstyle \footnotesize 
                 \begin{array}{c}}
\newcommand{\esb}{\end{array}}

\newtheorem{theorem}{{\sc Theorem}}[section]

 \newtheorem{lemma}{{\sc Lemma}}[section]

\usepackage{fancyhdr}
\usepackage[colorlinks,citecolor=blue]{hyperref}
\begin{document}
\title{The Existence of Dyon Solutions for Generalized Weinberg-Salam Model}
\author{Shouxin Chen\\School of Mathematics and Statistics\\Henan University\\
Kaifeng, Henan 475004, PR China\\ \\Yilu Xu\footnote{E-mail address:
104753180630@vip.henu.edu.cn.}\\School of Mathematics and Statistics\\Henan University\\Kaifeng, Henan 475004, PR China\\School of Economics\\Henan Kaifeng College of Science Technology and Communication\\
Kaifeng, Henan 475004, PR China}
\date{}
\maketitle

\textbf{Abstract}\\[-0.2cm]^^L
The generalized Weinberg-Salam model which is presented in a recent
study of Kimm, Yoon and Cho \cite{Kimm}, is arising in electroweak
theory. In this paper, we  prove the existence and asymptotic
behaviors  at infinity of static and radially symmetric dyon
solutions to the boundary-value problem of this model. Moreover, as
a by product, the qualitative properties of dyon solutions are also
obtained.
 The methods used here are the extremum principle,
  the Schauder fixed point
 theory and the shooting approach
  depending on one
shooting parameter. We provide an effective framework for
 constructing the dyon solutions in general dimensions and develop the existing
 results.
\begin{center}
\end{center}
 ~
  ~
\\~\\

\\
\mbox{}\hrule\mbox{}

\newpage
\tableofcontents

\baselineskip=15pt

\section{Introduction and main results}
This paper is concerned with the existence and asymptotic behaviors
at infinity of static and radially symmetric dyon solutions to the
boundary-value problem of the generalized Weinberg-Salam model (see
\cite{Kimm})
\begin{eqnarray}
\label{112} \mathcal{L}=-\left|\mathcal{D}_{\mu}\phi\right|^{2}
-\dfrac{\lambda}{2}\Big(\phi^{\dagger}\phi-\dfrac{\mu^{2}}{\lambda}\Big)^{2}
-\dfrac{1}{4}\vec{F}^{2}_{\mu\nu}-\dfrac{1}{4}G^{2}_{\mu\nu},
\end{eqnarray}
where $\mathcal{D}_{\mu}$ describes the covariant derivative of the
SU(2) subgroup only, $\dagger$ represents conjugate transpose,
$\vec{F}_{\mu\nu}$ and $G_{\mu\nu}$ are the gauge field strengths of
$SU(2)$ and $U(1)_{Y}$ with the potentials $\vec{A_{\mu}}$ and
$B_{\mu}$ and the corresponding coupling constants $g,\,g'$,
$\lambda,\,\mu$ are parameter, $\phi,\mathcal{D}_{\mu}\phi$ is
represented as
\begin{eqnarray}
&&\phi=\frac{1}{\sqrt{2}}\rho\xi,\,(\xi^{\dagger}\xi=1),\\
&&\mathcal{D}_{\mu}\phi=\Big(\partial_{\mu}-\mathrm{i}\frac{g}{2}\vec{\tau}\vec{A_{\mu}}
-\mathrm{i}\frac{g'}{2}B_{\mu}\Big)\phi
\triangleq\Big(D_{\mu}-\mathrm{i}\frac{g'}{2}B_{\mu}\Big)\phi,
\end{eqnarray}
where $\rho,\,\xi$ are the Higgs field and unit doublet,
respectively.

It is well known that theoretical physics and field theory, in
particular, provides a rich and challenging topic of study for
mathematics.~The study of these problems not only contributes to a
deeper understanding of physical concepts, theories, and the
relationships among them, but also provides new theories, methods,
and techniques for the development of mathematics.~In 1931, Dirac
has generalized the Maxwell's theory with his magnetic monopole
\cite{Dirac}.~Since then, magnetic monopoles have been the subject
of extensive research \cite{Lee}.~Dirac's monopole lives in the
classical Maxwell field theory for electromagnetism and carries
infinite energy \cite{Yang}.~By further studying Maxwell's equation,
Schwinger extended Dirac's idea of magnetic monopole \cite{tHooft},
and discovered and defined a new kind of particle like solution that
carries both electric charge and magnetic charge, which is called
dyon \cite{Schwinger}.~Dyon has a lot of important applications in
high-temperature superconductivity, quantum Hall effect and
superfluids \cite{Cho2}.~Generalized Yang-Mills theory has a
covariant derivative which contains both vector and scalar gauge
bosons \cite{Kimm}.~Based on this theory, some people construct an
SU(3) unified model \cite{tHooft} of weak and electromagnetic
interactions.~By using the NJL mechanism, the symmetry breaking can
be realized dynamically \cite{Hisano}.~The masses of W, Z
\cite{Barriola} are obtained and interactions between various
particles are the same as that of Weinberg-Salam(WS) model.~As for
monopoles, 't Hooft \cite{tHooft} and Polyakov \cite{Polyakov} have
shown the existence of finite-energy solutions in arrSO(3) Higgs
model.~However, the more relevant model for electromagnetic and weak
interactions is the SU(2) x U(1) model of Weinberg \cite{Weinberg}
and Salam.~The Weinberg-Salam electroweak model has important
theoretical significance and research value in classical field
theory \cite{Julia,Dokos,Manton,Rajaraman,Goddard}.

From a mathematical point of sight, the proof of the existence of a
magnetic monopole or dyon is a complex subject.~Their existences
rely either on explicit constructions in the self-dual limit
\cite{Actor,Bogomol,Prasad} or nonlinear functional analysis
\cite{Belavin,Tyupkin,Rawnsley,Jaffe,Maison}, as well as numerical
simulation \cite{Julia,Bais}.~Ever since Dirac \cite{Dirac}
predicted the existence of the monopole, the monopole has been an
obsession.~The Abelian monopole has been generalized to the
non-Abelian gauge theory by Wu and Yang \cite{Wu,Cho1} who showed a
non-Abelian monopole solution in the pure SU(2) gauge theory, and by
't Hooft \cite{tHooft} and Polyakov \cite{Polyakov} who have shown
that the SU(2) gauge theory allows a finite energy monopole solution
in Georgi-Glashow model as a topological soliton in the presence of
a triplet scalar source.~Moreover, the monopole in grand unification
has been constructed by Dokos and Tomaras \cite{Dokos}.~The
discovery of vortex solutions in the Weinberg-Salam model of the
electroweak interactions raises the possibility that such solutions
may exist in a wider class of field theories.~Indeed some time ago
Cho and Maison \cite{Cho2} have established that Weinberg-Salam
model and Georgi-Glashow model have exactly the same topological
structure, and demonstrated the existence of a new type of monopole
and dyon solutions in the standard Weinberg-Salam model.~Originally
the solutions of Cho and Maison were obtained by a numerical
integration \cite{Kimm}.~But a mathematically rigorous existence
proof has been established which endorses the numerical results.~And
the solutions are now referred to as Cho-Maison monopole and dyon
\cite{Schechter,YangBook,Yang}.~Up to now, many experimental results
have proved that Weinberg-Salam (WS) model \cite{Weinberg} is
correct in the current energy range.

In the spherically symmetric case, the
     two-point boundary value problem to \eqref{112}  is consisting
     of six nonlinear ordinary differential equations.~Namely,
     this is a  highly
      nonlinear and strongly coupled nonlinear ordinary differential
       system.
 Therefore it is very difficult to
        handle.~The purpose of our paper is to establish an existence
         theorem of the dyon solutions for the generalized Weinberg-Salam model \eqref{112}.
In fact, such a study was carried
  out in the earlier paper of Mcleod and the existence of the
   Wcinberg-Salam dyon was rigorously established by the method
   of calculus of variations in the article of Yang \cite{YangBook}.
   ~In a recent interesting work by Kimm, Yoon and Cho \cite{Kimm},
    three different ways is discussed to estimate the mass of the
    electroweak monopole and the differential equations governing
    the static radially symmetric the Cho-Maison monopole and dyon
     solutions are constructed in the generalized Weinberg-Salam model.

For the generators of $SU(2)$, we use the conventional Pauli
matrices $\tau^{\alpha}\,(\alpha=1,2,3)$. Then, $\eqref{112}$ is
reduced to the following equations of motion
\begin{equation}
\label{113} \left \{
\begin{array}{ll}
\partial^{2}\rho=\big|\mathcal{D}_{\mu}\xi\big|^{2}\rho+
\frac{\lambda}{2}\big(\rho^{2}-\frac{2\mu^{2}}{\lambda}\big)\rho,\\[3mm]
\mathcal{D}^{2}\xi=-2\frac{\partial_{\mu}\rho}{\rho}\mathcal{D}_{\mu}\xi
+\big[\xi^{\dagger}\mathcal{D}^{2}\xi+2\frac{\partial_{\mu}\rho}{\rho}
\big(\xi^{\dagger}\mathcal{D}_{\mu}\xi\big)\big]\xi,\\[3mm]
D_{\mu}\vec{F}_{\mu\nu}=\mathrm{i}\frac{g}{2}\rho^{2}\big[\xi^{\dagger}\vec{\tau}
\big(\mathcal{D}_{\nu}\xi\big)-\big(\mathcal{D}_{\nu}\xi\big)^{\dagger}
\vec{\tau}\xi\big],\\[3mm]
\partial_{\mu}G_{\mu\nu}=\mathrm{i}\frac{g'}{2}\rho^{2}\big[\xi^{\dagger}
\big(\mathcal{D}_{\nu}\xi\big)
-\big(\mathcal{D}_{\nu}\xi\big)^{\dagger}\xi\big].
\end{array}
\right.
\end{equation}
By the Abelian decomposition \cite{Kimm}
\begin{eqnarray}
\label{114}
\vec{\Phi}=\rho\hat{n},\,\,\,\hat{A}_{\mu}+\vec{W}_{\mu}=\vec{A}_{\mu},
\end{eqnarray}
we have
\begin{eqnarray}
\label{115}
&&\mathcal{L}=-\frac{1}{2}\big(\partial_{\mu}\rho\big)^{2}
-\frac{\rho^{2}}{2}\big|\hat{\mathcal{D}}_{\mu}\xi\big|^{2}
-\frac{\lambda}{8}\big(\rho^{2}-\rho_{0}^{2}\big)^{2}
-\frac{1}{4}\hat{F}_{\mu\nu}^{2} -\frac{1}{4}G_{\mu\nu}^{2}
-\frac{g}{2}\hat{F}_{\mu\nu}\cdot\big(\vec{W}_{\mu}\times\vec{W}_{\nu}\big)
\notag\\[2mm]
&&~~~~~~-\frac{1}{4}
\big(\hat{D}_{\mu}\vec{W}_{\nu}-\hat{D}_{\nu}\vec{W}_{\mu}\big)^{2}
-\frac{g^{2}}{8}\rho^{2}\big(\vec{W}_{\mu}\big)^{2}
-\frac{g^{2}}{4}\big(\vec{W}_{\mu}\times\vec{W}_{\nu}\big)^{2},
\end{eqnarray}
where
\begin{eqnarray*}
\label{116}
\hat{\mathcal{D}}_{\mu}=\partial_{\mu}-\mathrm{i}\frac{g}{2}\vec{\tau}\cdot\hat{A}_{\mu}
-\mathrm{i}\frac{g'}{2}B_{\mu}.
\end{eqnarray*}

To construct the desired solutions we enlarge $U(1)_{Y}$ and embed
it to another $SU(2)$.~And then, we introduce a hypercharged vector
field $X_{\mu}$ and a Higgs field $\sigma$, and generalize the
Lagrangian $\eqref{115}$ adding the following Lagrangian
\begin{eqnarray}
&&\Delta\mathcal{L}=-\frac{1}{2}\big|\tilde{D}_{\mu}X_{\nu}
-\tilde{D}_{\nu}X_{\mu}\big|^{2}
+\mathrm{i}g'G_{\mu\nu}X^{\ast}_{\mu}X_{\nu}
+\frac{1}{4}g'^{2}\big(X^{\ast}_{\mu}X_{\nu}-X^{\ast}_{\nu}X_{\mu}\big)^{2}
\notag\\[2mm]
&&~~~~~~~~~-\frac{1}{2}\left(\partial_{\mu}\sigma\right)^{2}
-g'^{2}\sigma^{2}\left|X_{\mu}\right|^{2}
-\frac{\kappa}{4}\Big(\sigma^{2}-\frac{m^{2}}{\kappa}\Big)^{2},
\end{eqnarray}
where \begin{eqnarray*}
\tilde{D}_{\mu}=\partial_{\mu}+\mathrm{i}g'B_{\mu}.
\end{eqnarray*}
Therefore, we can get the generalized Weinberg-Salam model in
electroweak theory as follows
\begin{eqnarray}
\label{000} &&\mathcal{L}=-\frac{1}{2}\big(\partial_{\mu}\rho\big)
-\frac{\rho^{2}}{2}\big|\hat{\mathcal{D}}_{\mu}\xi\big|^{2}
-\frac{\lambda}{8}\big(\rho^{2}-\rho_{0}^{2}\big)^{2}
-\frac{1}{4}\hat{F}_{\mu\nu}^{2} -\frac{1}{4}G_{\mu\nu}^{2}
-\frac{g}{2}\hat{F}_{\mu\nu}\cdot\big(\vec{W}_{\mu}\times\vec{W}_{\nu}\big)
-\frac{g^{2}}{8}\rho^{2}\big(\vec{W}_{\mu}\big)^{2}
\notag\\[2mm]
&&~~~~~~-\frac{1}{4}
\big(\hat{D}_{\mu}\vec{W}_{\nu}-\hat{D}_{\nu}\vec{W}_{\mu}\big)^{2}
-\frac{g^{2}}{4}\big(\vec{W}_{\mu}\times\vec{W}_{\nu}\big)^{2}
-\frac{1}{2}\big|\tilde{D}_{\mu}X_{\nu}
-\tilde{D}_{\nu}X_{\mu}\big|^{2}
+\mathrm{i}g'G_{\mu\nu}X^{\ast}_{\mu}X_{\nu}
\notag\\[2mm]
&&~~~~~~+\frac{1}{4}g'^{2}\big(X^{\ast}_{\mu}X_{\nu}-X^{\ast}_{\nu}X_{\mu}\big)^{2}
-\frac{1}{2}\left(\partial_{\mu}\sigma\right)^{2}
-g'^{2}\sigma^{2}\left|X_{\mu}\right|^{2}
-\frac{\kappa}{4}\Big(\sigma^{2}-\frac{m^{2}}{\kappa}\Big)^{2}.
\end{eqnarray}

In order to pursue a static radially symmetrie dyon solution, we
follow Kimm, Yoon and Cho to use the following general ansatz
\cite{Kimm}
\begin{eqnarray*}
\label{117}
&&\rho=\rho\left(r\right),\,\,\sigma=\sigma(r),\notag\\[2mm]
&&\xi=\mathrm{i}
\begin{pmatrix}
\sin\left(\frac{\theta}{2}\right)\mathrm{e}^{-\mathrm{i}\varphi}\\
-\cos\left(\frac{\theta}{2}\right)
\end{pmatrix},\notag\\[2mm]
&&X_{\mu}=\frac{i}{g'}\frac{h(r)}{\sqrt{2}}\mathrm{e}^{\mathrm{i}\varphi}
\left(\partial_{\mu}\theta+\mathrm{i}\sin\theta\partial_{\mu}\varphi\right),
\notag\\[2mm]
&&B_{\mu}=\frac{1}{g'}B(r)\partial_{\mu}t
-\frac{1}{g'}\left(1-\cos\theta\right)\partial_{\mu}\varphi,\notag\\[2mm]
&&\vec{A}_{\mu}=\frac{1}{g}A(r)\partial_{\mu}t\hat{r}
+\frac{1}{g}\left(f(r)-1\right)\hat{r}\times\partial_{\mu}\hat{r},
\end{eqnarray*}
where $\xi^{\dagger}\vec{\tau}\xi=-\hat{r}$, $(t,r,\theta,\varphi)$
are the spherically symmetric coordinates.

With the spherically symmetric ansatz the equations of motion for
functions $f(r),\rho(r),A(r)$, $B(r),h(r),\sigma(r),0<r<\infty$ are
reduced to
\begin{eqnarray}
\label{001} &f''-\dfrac{f^{2}-1}{r^{2}}f
=\left(\dfrac{g^{2}}{4}\rho^{2}-A^{2}\right)f,&\\[2mm]
\label{002} &\rho''+\dfrac{2}{r}\rho'-\dfrac{f^{2}}{2r^{2}}\rho
=-\dfrac{1}{4}(A-B)^{2}\rho+\dfrac{\lambda}{2}
\left(\rho^{2}-\dfrac{2\mu^{2}}{\lambda}\right)\rho,&\\[2mm]
\label{003} &A''+\dfrac{2}{r}A'-\dfrac{2f^{2}}{r^{2}}A
=\dfrac{1}{4}g^{2}\rho^{2}\left(A-B\right),&\\[2mm]
\label{004} &B''+\dfrac{2}{r}B'-\dfrac{2h^{2}}{r^{2}}B
=\dfrac{1}{4}g'^{2}\rho^{2}\left(B-A\right),&\\[2mm]
\label{005} &h''-\dfrac{h^{2}-1}{r^{2}}h
=\left(g'^{2}\sigma^{2}-B^{2}\right)h,&\\[2mm]
\label{006}
&\sigma''+\dfrac{2}{r}\sigma'-\dfrac{2h^{2}}{r^{2}}\sigma
=\kappa\left(\sigma^{2}-\dfrac{m^{2}}{\kappa}\right)\sigma.&
\end{eqnarray}

The boundary conditions for a regular field configuration can be
chosen as
\begin{eqnarray}
\label{007}
&f(0)=h(0)=1,\,A(0)=B(0)=\rho(0)=\sigma(0)=0,&\\[1mm]
\label{008} &f(\infty)=h(\infty)=0,\,A(\infty)=A_{0},
\,B(\infty)=B_{0},\,\rho(\infty)=\rho_{0},\,\sigma(\infty)=\sigma_{0},&
\end{eqnarray}
where $\sigma_{0}=\sqrt{\frac{m^{2}}{\kappa}}$,
$\rho_{0}=\mu\sqrt{\frac{2}{\lambda}}$,
$\lambda,\,\mu,\,g,\,g',\,\kappa,\,m$ are parameters,
$A_{0},$\,$B_{0}$ are given positive constants.

For the above nonlinear ordinary differential equation with
two-point boundary value problem, inspired by the literature
\cite{McLeod1,Bizon,Hastings,Smoller,Wang,Mcleod2}, we develop the
methods and techniques in which we prove the existence of the
solution and study the related properties of the solution.

To do so, we require parameters to satisfy the following two
assumptions
\begin{eqnarray*}
&&{\rm(H1)}\,\frac{1}{4}g\rho^{2}_{0}>A^{2}_{0},\,
g'^{2}\sigma^{2}_{0}>B^{2}_{0},\\[2mm]
&&{\rm(H2)}\,B_{0}=A_{0}.
\end{eqnarray*}
Then the main results of this paper are stated as follows.
\begin{theorem}[Existence of solutions to the boundary-value problem
]\label{th1} Under the assumption (H1) and (H2), the Weinberg-Salam
dyon equations $\eqref{001}-\eqref{006}$ have a family of finite
energy smooth solutions which satisfy the radial symmetry
properties.~The obtained solution configuration functions
$(f,\rho,A,B,h,\sigma)$ have the properties that

{\rm(1)} The function $f\,,A\,,B\,,h\in
C^{1}\left([0,+\infty)\right)$ and $f'(0)=h'(0)=0${\rm;}

{\rm(2)} $0\leqslant
f(r),h(r)\leqslant1,$\,$\rho^{2}(r)\leqslant\rho_{0}^{2}
+\frac{A_{0}^{2}}{2\lambda},$\,$A(r)\leqslant A_{0},B(r)\leqslant
A_{0},$\,$B(r)\geqslant A(r)$ for all $r\geqslant0{\rm;}$

{\rm(3)} $r\rho(r)\,,rA(r)\,,rB(r)\,,r\sigma(r)$ are increasing,
$f(r)\,,h(r)\,,r^{-1}B(r)\,,r^{-1}A(r)$ are decreasing{\rm;}

{\rm(4)} $r^{-k}\rho(r)$ is decreasing as long as
$\rho\leqslant\rho_{0}$, where $k=\frac{1}{2}(\sqrt{3}-1);$
$r^{-2}\sigma(r)$ is decreasing as long as
$\sigma\leqslant\sigma_{0}${\rm;}

{\rm(5)} $r^{-1}B(r),$\,$r^{-1}A(r)$ is bounded as $r\rightarrow0$.
\end{theorem}

About the asymptotics of the solutions as \ $r\to\infty$, we have
the following result.
\begin{theorem}[The
asymptotic exponential decay]\label{th2} As $r\rightarrow \infty$,
there hold the sharp asymptotic estimates
\begin{eqnarray}
&f(r)=O(e^{-\kappa(1-\varepsilon)r}),&\\\label{01}
&\rho(r)=\rho_{0}+O(r^{-1}e^{-\sqrt{2}\mu_{0}(1-\varepsilon)r}),&\\\label{02}
&A(r)=A_{0}+O(r^{-1}),& \\\label{03}
&h(r)=O(e^{-\zeta(1-\varepsilon)r}),&\\\label{04}
&B(r)-A(r)=O(r^{-1}e^{-\nu_{0}(1-\varepsilon)r}),&\\\label{05}
&\sigma(r)=\sigma_{0}+O(r^{-1}e^{-\sqrt{2}\xi(1-\varepsilon)r}),&\label{06}
\end{eqnarray}
where $\varepsilon>0$ can be taken to be arbitrarily small,
$\rho_{0}=\mu\sqrt{\frac{2}{\lambda}}$,
$\nu=\dfrac{\rho_{0}}{2}\sqrt{g^{2}+g'^{2}}$ and the decay exponents
are defined by the expressions
\begin{eqnarray*}
  &\kappa=\sqrt{\frac{1}{4}g^{2}\rho_{0}^{2}-A_{0}^{2}},\,\,
  \mu_{0}=\min\{\mu,\sqrt{2}\kappa,\frac{\nu}{\sqrt{2}}\},&\\[2mm]
  &\zeta=\sqrt{g'^{2}\sigma_{0}^{2}-A_{0}^{2}},\,\,
  \nu_{0}=\min\{2\kappa,\nu\},\,\,
  \xi=\sqrt{\kappa}\sigma_{0}.&
\end{eqnarray*}

\noindent\textbf{Remark:} In the proof of Theorem \ref{th1}, we will
give asymptotic estimates for the solution of the above problem when
$r\rightarrow0$, and we omit them here.

\end{theorem}

As the end of this section, we state the arrangement of this paper
as follows.~In Section 2, we first give  a series of lemmas as the
primary works for proving our main results. ~And then, by using the
shooting method and the Sturm comparison principle,  we prove the
existence of the dyon solution of each second-order nonlinear
ordinary differential equation and establish the qualitative
properties of  solutions.~In Section 3, by using the Schauder fixed
point theory, we study the two-point boundary value problem
\eqref{001}-\eqref{008} and obtain the asymptotic behaviors of the
solutions at infinity.

\section{Primaries}

The proof of the theorem depends on a fixed-point argument, and we
outline this in the statement of a series of lemmas, which will  be
proved in sequence.
\begin{lemma}\label{le1}
Given a pair of fuctions $(\rho(r),A(r),h(r))\in C\left([0,+\infty)\right)$ such that

{\rm(1)} $h'(0)=h(\infty)=0$, $h(0)=1$, $\rho(\infty)=\rho_{0}$, $A(\infty)=A_{0};$

{\rm(2)} When $r\rightarrow0$, $r^{-k}\rho(r),\,\,r^{-1}A(r),\,\,r^{-1}(h(r)-1)$ exist finite limits{\rm;}

{\rm(3)} $h(r)$ is decreasing, $r\rho(r),\,\,rA(r)$ are increasing{\rm;}

{\rm(4)} $r^{-k}\rho(r),\,\,r^{-1}A(r)$ is decreasing so long as $r^{\alpha}\leqslant\frac{\rho_{0}}{R^{*}}\leqslant1$, that is $\rho(r)\leqslant\rho_{0}${\rm;}

{\rm(5)} $\rho^{2}(r)\leqslant\rho_{0}^{2}+\frac{A_{0}^{2}}{2\lambda}$, $A(r)\leqslant A_{0}${\rm;}

{\rm(6)}\ $r^{-\alpha}\rho(r),\,r^{-\alpha}A(r),\,r^{-\alpha}(h(r)-1)\leqslant R^{*},\,\forall r\leqslant1${\rm;}

\noindent where $0<\alpha<k=\frac{1}{2}(\sqrt{3}-1)$, $R^{*}=R^{*}\left(\rho_{0},\lambda,\mu,A_{0},g,g'\right)$ is a positive constant $($without loss of generality, we could assume $R^{*}\geqslant\rho_{0})$, then we can find a unique continuously differentiable function f satisfying equation $\eqref{001}$ and the conditions as follow

{\rm(1)} $f(0)=1,\,f(\infty)=0,\,f'(0)=0$, $f(r)$ is decreasing{\rm;}

{\rm(2)} $|r^{-2}\psi(r)|\leqslant N\left({R^{*}}^{\frac{2}{1+\alpha}}+r^{2\alpha}{R^{*}}^{2}\right)$\,$(\forall r\leqslant1)$, where $\psi(r)=f(r)-1$, $N=N(\rho_{0},\lambda,\mu,A_{0},g,g')$ is a positive constant.
\end{lemma}

\begin{lemma}\label{le2}
Given the functions $(\rho(r),A(r),h(r))$ as in Lemma \ref{le1}, and the associated function $f(r)$, then we can find a unique function $B(r)\in C^{1}\left([0,+\infty)\right)$ satisfying the equation $\eqref{004}$ and the conditions
\begin{equation*}
B(0)=0,\,\,B(\infty)=B_{0},\,\,A(r)\leqslant B(r)\leqslant A_{0},
\,\,|r^{-1}B(r)|\leqslant R^{*}\,\,(\forall r\leqslant1).
\end{equation*}
\end{lemma}

\begin{lemma}\label{le3}
Given the functions $(\rho(r),A(r),h(r))$ as in Lemma \ref{le1}, and the associated function $f(r)$, then we can find a unique function $\sigma(r)\in C^{1}\left([0,+\infty)\right)$ satisfying the equation $\eqref{006}$ and the conditions
\begin{eqnarray*}
&\sigma(0)=0,\,\,\sigma(\infty)=\sigma_{0},\,\,r\sigma(r)\,\,is\,\,increasing,&\\
&r^{-1}\sigma(r)\,\,is\,\,decreasing\,\,so\,\,long\,\,as\,\,\sigma\leqslant\sigma_{0},
\,\,|r^{-1}\sigma(r)|\leqslant R^{*}\,(\forall r\leqslant1).&
\end{eqnarray*}
\end{lemma}

\begin{lemma}\label{le4}
Given the function $(\rho(r),A(r),h(r))$ as in Lemma \ref{le1}, and the associated function $B(r)$,\,$\sigma(r)$ from Lemmas \ref{le2},\,\ref{le3}, then we can find a unique function $\tilde{h}(r)\in C^{1}\left([0,+\infty)\right)$ satisfying the equation
\begin{equation}\label{0007}
\tilde{h}''-\dfrac{\tilde{h}^{2}-1}{r^{2}}\tilde{h}
=\left(g'^{2}\sigma^{2}-B^{2}\right)\tilde{h}
\end{equation}
and the conditions $\tilde{h}(r)$ is decreasing,
\begin{eqnarray*}
\tilde{h}(0)=0,\,\,\tilde{h}(\infty)=1,\,\,\tilde{h}'(0)=0,\,\,
|r^{-2}(\tilde{h}(r)-1)|\leqslant N({R^{*}}^{\frac{2}{1+\alpha}}+r^{2\alpha}{R^{*}}^{2})
\,(\forall r\leqslant1).
\end{eqnarray*}
\end{lemma}

\begin{lemma}\label{le5}
Given the function $(\rho(r),A(r),h(r))$ as in Lemma \ref{le1}, and the associated function $f(r)$,\,$B(r)$ from Lemmas \ref{le1},\,\ref{le2}, then we can find a unique function $\tilde{\rho}(r)\in C^{1}\left([0,+\infty)\right)$ satisfying the equation
\begin{equation}
\label{0009}
\left(r\tilde{\rho}\right)''-\dfrac{f^{2}}{2r}\tilde{\rho}
=-\dfrac{1}{4}r\left(A-B\right)^{2}\tilde{\rho}+\dfrac{\lambda}{2}r\tilde{\rho}
\left(\tilde{\rho}^{2}-\dfrac{2\mu^{2}}{\lambda}\right)
\end{equation}
and the conditions
\begin{eqnarray*}
&\tilde{\rho}(0)=0,\,\,\tilde{\rho}(\infty)=\rho_{0},
\,\,\tilde{\rho}^{2}(r)\leqslant\rho_{0}^{2}+\frac{A_{0}^{2}}{2\lambda},\,\,
r\tilde{\rho}(r)\,\,is\,\,increasing,&\\
&r^{-k}\tilde{\rho}(r)\,\,is\,\,decreasing\,\,so\,\,long\,\,as\,\,
\tilde{\rho}(r)\leqslant\rho_{0},\,\,
|r^{-k}\tilde{\rho}(r)|\leqslant R_{2}^{*}\,(\forall r\leqslant1),\,\,where\,\,k=\frac{1}{2}(\sqrt{3}-1).&
\end{eqnarray*}
\end{lemma}

\begin{lemma}\label{le6}
Given the function $(\rho(r),A(r),h(r))$ as in Lemma \ref{le1}, and the associated function $f(r)$,\,$B(r)$ from Lemmas \ref{le1},\,\ref{le2}, then we can find a unique function $\tilde{A}(r)\in C^{1}\left([0,+\infty)\right)$ satisfying the equation
\begin{equation}\label{0011}
\left(r\tilde{A}\right)''-\dfrac{2}{r}f^{2}\tilde{A}
=-\dfrac{1}{4}g^{2}\tilde{\rho}^{2}r\left(B-\tilde{A}\right)
\end{equation}
and the conditions
\begin{eqnarray*}
&\tilde{A}(0)=0,\,\,\tilde{A}(\infty)=A_{0},\,\,
r\tilde{A}(r)\,\,is\,\,increasing,&\\
&r^{-1}\tilde{A}(r)\,\,is\,\,decreasing,
\,\,\tilde{A}(r)\leqslant B(r),
\,\,r^{-1}\tilde{A}(r)\leqslant R_{3}^{*}
\,(\forall r\leqslant1).&
\end{eqnarray*}
\end{lemma}

The final theorem is then just a matter of constructing a mapping and showing it has a fixed point, and this is proved in the final section and so does the asymptotics of the solutions.

\subsection{Proof of Lemma \ref{le1} (Existence and uniqueness of\ $f(r)$)}

In this subsection, we prove in two steps.~To solve the two-point
boundary value problem, we use a single parameter shooting
method.~When we do this, we need to consider the initial value
problem
\begin{eqnarray}
\label{0001}
f''-\dfrac{f^{2}-1}{r^{2}}f
=\left(\frac{g^{2}}{4}\rho^{2}-A^{2}\right)f,\,\,f(0)=1.
\end{eqnarray}

Firstly, we prove the existence of local solutions to the initial
value problem $\eqref{0001}$.

\begin{lemma}\label{le7}
There exists a locally continuous solution of the initial value problem $\eqref{0001}$ near $x=0$.
\end{lemma}
{\bf Proof } Using the idea from the literature \cite{Bizon} and \cite{Smoller}, this lemma will be proved. Firstly, we introduce a new variable $w=f'$, and rewrite the equation $\eqref{0001}$ as the first-order ordinary differential equation system
\begin{eqnarray}
\label{200}
&&f'=w,\\[1mm]
&&\left(r^{2}w\right)'=2rw+\left[\left(f^{2}-1\right)
+r\left(\frac{g^{2}}{4}\rho^{2}-A^{2}\right)\right]f.
\end{eqnarray}
Then consider the space $\mathcal{X}$ as follows
\begin{eqnarray*}
\mathcal{X}=\left\{(f(r),w(r))\in C\left([0,x]\right)|\,\left\|f(r)-1\right\|\leqslant1,
\,\left\|w(r)\right\|\leqslant1,\,
\forall r\in[0,x]\right\}
\end{eqnarray*}
where $\left\|h\right\|=\sup\left\{|h(x)|:0\leqslant r\leqslant x\right\}$.~It is clear that $\mathcal{X}$ is a complete normed linear space if we take as metric the maximum of the two components.~We define a map \,$T{\rm:}\,(f,w)\rightarrow(T_{1},T_{2})$ on $\mathcal{X}$ where
\begin{eqnarray}
\label{201}
&&T_{1}=1+\int_0^r wds,\\[1mm]
&&T_{2}=\frac{1}{r^{2}}\int_0^r 2sw+\left[\left(f^{2}-1\right)
+s\left(\frac{g^{2}}{4}\rho^{2}-A^{2}\right)\right]fds.
\end{eqnarray}

One verifies easily that $T$ does in fact take $\mathcal{X}$ to $\mathcal{X}$ and that $T$ is a contracting map if $r$ is sufficiently small, and that a fixed point of $T$ is a solution to our equation.~Therefore the lemma follows.
\hfill$\Box$\vskip7pt

Setting $\psi=f-1$, then we may rewrite equation $\eqref{0001}$ as
\begin{eqnarray}
\label{011}
\psi''=\frac{2\psi}{r^{2}}
+\left(\frac{1}{4}g^{2}\rho^{2}-A^{2}\right)(\psi+1)
+\frac{\psi^{3}+3\psi^{2}}{r^{2}}.
\end{eqnarray}
Using the basic theory of ordinary differential equations and initial value condition, the differential equation can be transformed into the following integral equation
\begin{eqnarray}
\label{012}
\psi(r)=Cr^{2}+\frac{1}{3}\int_0^r\left(r^{2}s^{-1}-r^{-1}s^{2}\right)
\left\{\left(\frac{1}{4}g^{2}\rho^{2}-A^{2}\right)(\psi+1)
+\frac{3\psi^{2}+\psi^{3}}{s^{2}}\right\}ds,
\end{eqnarray}
where $C$ is an arbitrary constant.~Equation can be solved by lemma \ref{le7}, at least
for r sufficiently small. And according to $\eqref{012}$, we obtain a solution with
\begin{equation}
\psi(r)=Cr^{2}+O\left(r^{2+2\alpha}\right)\,(r\rightarrow0^{+}).
\end{equation}

Applying the extension theorem of the solution of the ordinary
differential equation, the solution can be extended to its maximum
existence interval $[0,R_{C})$, where either $R_{C}=\infty$ or
$\lim\limits_{r\rightarrow R_{C}^{-}}\psi(r)=\infty$.~By the
continuous dependence of the solution on the parameters theorem we
obtain that the solution $\psi$ depends continuously on the
parameter $C$.

Next we show the existence of the global solution to the boundary
value problem, then we consider $C>0$ and define the sets
$\mathcal{S}_{1}^{f},\mathcal{S}_{2}^{f},\mathcal{S}_{3}^{f}$ as
follows
\begin{eqnarray*}
\mathcal{S}_{1}^{f}&=&\left\{C<0:\,\psi'(r;C)\,\,\mbox{becomes positive before}\,\,\psi(r;C)\,\,\mbox{reaches}\,-1\right\},\\[1mm]
\mathcal{S}_{2}^{f}&=&\left\{C<0:\,\psi(r;C)\,\,\mbox{crosses}\,-1\,
\,\mbox{before}\,\,\psi'(r;C)\,\,\mbox{becomes}\,\,0\right\},\\[1mm]
\mathcal{S}_{3}^{f}&=&\left\{C<0:\,\forall r>0,\psi'(r;C)\leqslant0,-1<\psi(r;C)<0\right\}.
\end{eqnarray*}
Obviously, it follows from the construction of the set
\begin{equation*}
\mathcal{S}_{1}^{f}\cup\mathcal{S}_{2}^{f}\cup\mathcal{S}_{3}^{f}=\mathcal{S}^{f},\,
\,\,\mathcal{S}_{1}^{f}\cap\mathcal{S}_{2}^{f}
=\mathcal{S}_{2}^{f}\cap\mathcal{S}_{3}^{f}
=\mathcal{S}_{3}^{f}\cap\mathcal{S}_{1}^{f}=\emptyset.
\end{equation*}

\begin{lemma}\label{le11}
The sets $\mathcal{S}_{1}^{f},\mathcal{S}_{2}^{f}$ are both open and
nonempty.
\end{lemma}
{\bf Proof } Firstly, $\mathcal{S}_{1}^{f}$ contains $C$
small.~Inserting $C=0$ into the equation $\eqref{012}$, we obtain
the equations as follows
\begin{eqnarray}
\label{018}
\psi(r)=\frac{1}{3}\int_0^r(r^{2}s^{-1}-r^{-1}s^{2})
\left\{\left(\frac{1}{4}g^{2}\rho^{2}-A^{2}\right)(\psi+1)
+\frac{3\psi^{2}+\psi^{3}}{s^{2}}\right\}ds,
\end{eqnarray}
\begin{eqnarray}
\label{019}
\psi'(r)=\frac{1}{3}\int_0^r(2rs^{-1}+r^{-2}s^{2})
\bigg \{\bigg(\frac{1}{4}g^{2}\rho^{2}-A^{2}\bigg)(\psi+1)
+\frac{3\psi^{2}+\psi^{3}}{s^{2}}\bigg \}ds.
\end{eqnarray}

When $r>0$ is sufficiently small, it is obvious that
$\psi'>0$.~According to $\eqref{0001}$, we obtain that
$\frac{1}{4}g\rho^{2}-A^{2}>0$.~In addition,
$\psi+1>0,\,\frac{3\psi^{2}+\psi^{3}}{s^{2}}>0$ are bounded.~Hence
$\psi>0$ when $r>0$ is sufficiently small.~By the continuous
dependence of $\psi$ and $\psi'$ on the parameter $C$, we obtain
that $\psi>0,\,\psi'>0$ when $C<0,\,r>0$ are both sufficiently
small.~Therefore $C\in\mathcal{S}_{1}^{f}$ and the nonemptyness of
$\mathcal{S}_{1}^{f}$ is established.~It is evident that
$\mathcal{S}_{1}^{f}$ is open because the continuous dependence of
$\psi$ and $\psi'$ on the parameter $C$.

Secondly, $\mathcal{S}_{2}^{f}$ contains $C$ large.~Here we
introduce a transformation to consider a modified variable
$t=|C|^{\frac{1}{2}}r$ in $\eqref{012}$ so that it becomes
\begin{eqnarray}
\label{020}
\psi(t)=-t^{2}+\frac{1}{3}\int_0^t(t^{2}\tau^{-1}-t^{-1}\tau^{2})\bigg \{|C|^{-1}\bigg(\frac{1}{4}g^{2}\rho^{2}-A^{2}\bigg)(\psi+1)
+\frac{3\psi^{2}+\psi^{3}}{\tau^{2}}\bigg\}d\tau,
\end{eqnarray}
\begin{eqnarray}
\label{021}
\psi'(t)=-2t+\frac{1}{3}\int_0^t(2t\tau^{-1}+t^{-2}\tau^{2})
\bigg\{|C|^{-1}\bigg(\frac{1}{4}g^{2}\rho^{2}-A^{2}\bigg)(\psi+1)
+\frac{3\psi^{2}+\psi^{3}}{\tau^{2}}\bigg\}d\tau.
\end{eqnarray}
The differential equation corresponding to $\eqref{020}$ is
\begin{eqnarray}
\label{022}
\psi_{tt}-\frac{\psi(\psi+1)(\psi+2)}{t^{2}}
=\bigg(\frac{1}{4}g^{2}\rho^{2}-A^{2}\bigg)(\psi+1)|C|^{-1}.
\end{eqnarray}
Since $r^{-\alpha}\rho,r^{-\alpha}A\leqslant R^{*},\forall
r\leqslant1$, we obtain that
\begin{eqnarray*}
  \rho\leqslant|C|^{\frac{-\alpha}{2}}t^{\alpha}R^{*},\,
  A\leqslant|C|^{\frac{-\alpha}{2}}t^{\alpha}R^{*},\,\forall t\leqslant|C|^{\frac{1}{2}}.
\end{eqnarray*}
From the continuity of the function $\psi(r)$ on $[0,R]$, then for all $t\leqslant|C|^{\frac{1}{2}}$, there is an $M>0$ such that $|\psi(t)|\leqslant M$.~Hence we obtain the estimate
\begin{eqnarray*}
|C|^{-1}\bigg|\frac{1}{4}g^{2}\rho^{2}-A^{2}\bigg|(\psi(t)+1)
\leqslant|C|^{-1-\alpha}\tau^{2\alpha}R^{*}M\bigg(\frac{1}{4}g^{2}+1\bigg).
\end{eqnarray*}
According to
\begin{eqnarray*}
\bigg(\dfrac{1}{4}g^{2}\rho^{2}-A^{2}\bigg)(\psi+1)|C|^{-1}\rightarrow0\,\,as
\,\,C\rightarrow-\infty,
\end{eqnarray*}
we obtain that $\eqref{020}$ is equivalent to the differential
equation
\begin{eqnarray}
\label{023}
\frac{d^{2}\psi}{dt^{2}}=\frac{\psi(\psi+1)(\psi+2)}{t^{2}},\,\,\forall t\in[0,R],
\end{eqnarray}
where $R$ is a positive constant.~It is clear that the solution of $\eqref{022}$ as follows
\begin{eqnarray*}
\psi\sim-t^{2},\,\,(r\rightarrow0,\,C\rightarrow-\infty).
\end{eqnarray*}

Thus the solution $\psi$ could cross $-1$ at $t_{0}=1+\varepsilon_{0}(\varepsilon_{0}>0)$.~Noting that $\psi_{t}(t)<0,\,\forall t\in(0,t_{0}]$, then $\psi(r;C)$ crosses $-1$ before $\psi'(r;C)$ becomes $0$.~Therefore there exists some positive constant $N=N(\rho_{0},\lambda,\mu,A_{0},g,g')$ independent of the choice of $\rho$, $A$, and $R^{*}$ such that $C\in\mathcal{S}_{2}^{f}$ when $R^{*}|C|^{-1-\alpha}<\frac{1}{N}$.~In other words, $\mathcal{S}_{2}^{f}$ is nonempty.~The fact that $\mathcal{S}_{2}^{f}$ is open is self-evident.
\hfill$\Box$\vskip7pt

Since the connected set $C<0$ cannot consist of two open disjoint non-empty sets, there must be some value of $C$ in neither $\mathcal{S}_{1}^{f}$ nor $\mathcal{S}_{2}^{f}$.~For this value of $C$, say $C_{0}$, we have a solution with $\psi'(r;C)\leqslant0,\,-1<\psi(r;C)<0$.

\begin{lemma}\label{le12}
The solution corresponding to the parameter $C_{0}$ in $\mathcal{S}_{3}^{f}$ satisfies $\lim\limits_{r\rightarrow\infty}f(r)=0$.
\end{lemma}
{\bf Proof } Since $-1<\psi(r;C_{0})<0$ and
$\psi'(r;C_{0})\leqslant0$, $\forall r>0$, we obtain that
$\lim\limits_{r\rightarrow\infty}f(r;C_{0})\triangleq L\geqslant0$.
It is obvious that $f(r)=O(e^{-\sqrt{\kappa}r})(r\rightarrow\infty)$
because $f''\sim(\frac{g^{2}}{4}\rho^{2}-A^{2})f$ as
$r\rightarrow\infty$.~Noting that
\begin{eqnarray*}
\lim\limits_{r\rightarrow\infty}\kappa(r)
=\lim\limits_{r\rightarrow\infty}\left(\frac{g^{2}}{4}\rho^{2}-A^{2}\right)
=\frac{g^{2}}{4}\rho_{0}^{2}-A_{0}^{2}>0,
\end{eqnarray*}
we have $L=\lim\limits_{r\rightarrow\infty}f(r)=0$.
\hfill$\Box$\vskip7pt

Finally, we want to prove that the solution\ $f$ is the only solution satisfying the
conditions that
\begin{equation*}
f'(0)=0,\,\,\,|r^{-2}\psi(r)|\leqslant N\left({R^{*}}^{\frac{2}{1+\alpha}}+r^{2\alpha}{R^{*}}^{2}\right)
\,(\forall r\leqslant1).
\end{equation*}

\begin{lemma}\label{le13}
The solution for the given parameter $C_{0}$ is unique.
\end{lemma}
{\bf Proof } Suppose otherwise that there are two solutions $f_{1},f_{2}$, and set $\Psi(r)=f_{2}(r)-f_{1}(r)$. Then the function $\Psi(r)$ satisfies the boundary condition $\Psi(0)=\Psi(\infty)=0$ and the equation
\begin{eqnarray}
\label{027}
r^{2}\Psi''(r)&=&\left[(f_{1}^{2}(r)-1)
+\left(\frac{1}{4}g^{2}\rho^{2}(r)-A^{2}(r)\right)r^{2}
+(f_{2}^{2}(r)+f_{1}(r)f_{2}(r))\right]\Psi(r)\notag\\[1mm]
&\triangleq& Q(r)\Psi(r),\,\,\,0<r<+\infty.
\end{eqnarray}

Without loss of generality, we assume $\Psi>0$ when $r>0$ is
sufficiently small.~According to $\eqref{0001}$, we obtain that
\begin{eqnarray}
\label{026}
  r^{2}f_{i}''(r)&=&\left[(f_{i}^{2}(r)-1)
  +\left(\frac{1}{4}g^{2}\rho^{2}(r)-A^{2}(r)\right)r^{2}\right]f_{i}(r)\notag\\[1mm]
  &\triangleq& q_{i}(r)f_{1}(r),\,\,0<r<+\infty,\,\,(i=1,2).
\end{eqnarray}

Since $Q(r)-q_{1}(r)>0$, Applying the Sturm-Picone comparison
theorem to $\eqref{027}-\eqref{026}$, we conclude that $f_{1}$ have
more zero points than $\Psi$.~Note that $f_{1}(r)\neq0$ as
$r\in[0,+\infty)$, then we have $\Psi(r)\neq0$ for all
$r\in[0,+\infty)$, which contradicts $\Psi(0)=0$.
\hfill$\Box$\vskip7pt

\begin{lemma}\label{le14}
If $C_{0}\in\mathcal{S}_{3}^{f}$, for $r\leqslant1$, we can find $N=N(\rho_{0},\lambda,\mu,A_{0},g,g')>0$ such that
\begin{eqnarray*}
|r^{-2}\psi(r)|\leqslant N\left({R^{*}}^{\frac{2}{1+\alpha}}+r^{2\alpha}{R^{*}}^{2}\right),
\end{eqnarray*}
where $\psi=f-1$.~Furthermore, we have $f'(0)=0$.
\end{lemma}
{\bf Proof } According to $\eqref{012}$, $A(s)\leqslant s^{\alpha}R^{*},\,\,\rho(s)\leqslant s^{\alpha}R^{*},\,\,\forall s\leqslant1,$ and $|C_{0}|^{1+\alpha}\leqslant N{R^{*}}^{2}$ $($because $C_{0}\in\mathcal{S}_{3}^{f})$, we arrive at
\begin{eqnarray}
\label{063}
\left|\frac{\psi}{r^{2}}\right|
&\leqslant&|C_{0}|+\frac{1}{3}\int_0^r\left|s^{-1}
\left(\frac{1}{4}g^{2}+1\right)s^{2\alpha}{R^{*}}^{2}\right|ds\notag\\[1mm]
&\leqslant&\left(N{R^{*}}^{2}\right)^{\frac{1}{1+\alpha}}
+\frac{(\frac{1}{4}g^{2}+1)}{3}\int_0^r s^{2\alpha-1}{R^{*}}^{2}ds\notag\\[2mm]
&\leqslant&N'\left({R^{*}}^{\frac{2}{1+\alpha}}+r^{2\alpha}{R^{*}}^{2}\right)
\leqslant N''\left({R^{*}}^{\frac{2}{1+\alpha}}+{R^{*}}^{2}\right)\triangleq{R^{*}},
\,\,\forall r\leqslant1,
\end{eqnarray}
where $N'=\max\{N^\frac{1}{1+\alpha},
\frac{(\frac{1}{4}g^{2}+1)}{6\alpha}\}$.~Using the above result, we
have $\psi(r)=O(r^{2})(r\rightarrow0^{+})$.~By the defination of the
derivative of the function $\psi(r)$ at $r=0$, we obtain that
\begin{equation*}
  \psi'(0)=\lim\limits_{r\rightarrow0^{+}}\frac{\psi(r)-\psi(0)}{r}
  =\lim\limits_{r\rightarrow0^{+}}\frac{O(r^{2})}{r}=0.
\end{equation*}

Namely, $f'(0)=0$.
This completes the proof of the Lemma \ref{le1}.
\hfill$\Box$\vskip7pt

\subsection{Proof of Lemma \ref{le2} (Existence and uniqueness of\ $B(r)$)}

By noting that we can rewrite $\eqref{004}$ as
\begin{eqnarray}
\label{030}
r^{2}B''+2rB'-2B=2(h^{2}-1)B+\dfrac{1}{4}g'^{2}\rho^{2}r^{2}(B-A).
\end{eqnarray}

Similarly, differential equation $\eqref{030}$ can be transformed into the integral equation form
\begin{eqnarray}
\label{031}
B(r)=br+\frac{1}{3}\int_0^r\left(r-r^{-2}s^{3}\right)\left [2(h^{2}-1)B+\dfrac{1}{4}g'^{2}\rho^{2}s^{2}(B-A)\right]ds,
\end{eqnarray}
where $b$ is an arbitrary constant, which can be solved by Picard iteration to give a locally continuous solution of the initial value problem that exists at least for $r$ sufficiently small, and we obtain a solution with $B(r)=br+O(r^{2+2\alpha})$, and the solution continuously depends on the parameter $b$.

Same as previous section, we introduce three sets
\begin{eqnarray*}
\mathcal{S}^{B}_{1}&=&\left\{b>0:\,B(r;b)<A\,\,\mbox{before}\,\,B(r;b)
\,\,\mbox{reaches}\,\,A_{0}\right\},\\[1mm]
\mathcal{S}^{B}_{2}&=&\left\{b>0:\,B(r;b)\,\,\mbox{crosses}\,\,A_{0}\,\,\mbox{before}
\,\,B(r;b)=A\,\right\},\\[1mm]
\mathcal{S}^{B}_{3}&=&\left\{b>0:\,\forall r>0,B(r;b)\leqslant A_{0},B(r;b)\geqslant A\right\}.
\end{eqnarray*}

Obviously,
$\mathcal{S}^{B}_{1}\cup\mathcal{S}^{B}_{2}\cup\mathcal{S}^{B}_{3}=\mathcal{S}^{B},\,
\mathcal{S}^{B}_{1}\cap\mathcal{S}^{B}_{2}=\mathcal{S}^{B}_{2}\cap\mathcal{S}^{B}_{3}
=\mathcal{S}^{B}_{3}\cap\mathcal{S}^{B}_{1}=\emptyset$.

\begin{lemma}\label{le21}
The set $\mathcal{S}_{1}^{B},\mathcal{S}_{2}^{B}$ are both open and nonempty.
\end{lemma}
{\bf Proof } According to the integral equation $\eqref{031}$, if\ $b=0$, since
\begin{eqnarray}
\label{034}
B'(r)=\frac{1}{3}\int_0^r \left(1+2r^{-3}s^{3}\right)\left [2(h^{2}-1)B+\dfrac{1}{4}g'^{2}\rho^{2}s^{2}(B-A)\right]ds,
\end{eqnarray}
by iteration we obtain that $B'<0$ when $r>0$ is sufficiently
small.~If $r>0$ is near zero, it is clear that $B-A<0$ because
$0<A<1$ and $ B(0)=0$.~By the continuous dependence of $B$ on the
parameter $b$ we obtain that $B-A<0$ for the sufficiently small
$b,r>0$.~Let $r=r_{0}$ be the first value satisfying $B(r)=A(r)$,
then
\begin{eqnarray*}
  B'(r)\bigg|_{r=r_{0}}
  =\left(b+\frac{1}{3}\int_0^r\left[
  \left(1+2r^{-3}s^{3}\right)\left(h^{2}-1\right)A\right]ds\right)\bigg|_{r=r_{0}}>0.
\end{eqnarray*}

Obviously,
\begin{eqnarray*}
0<b<\left|\frac{1}{3}\left\{\int_0^r
\left[\left(1+2r^{-3}s^{3}\right)\left(h^{2}-1\right)A\right]ds\right\}
\bigg|_{r=r_{0}}\right|,
\end{eqnarray*}
for all $b\in(0,A_{0}r_{0})$, then $B'|_{r=r_{0}}<0$, which contradicts $B'|_{r=r_{0}}>0$.~Hence, $B<A$ when $0<b\leqslant A_{0}r_{0}$.~The following proves that the above $b$ can be taken.~Since
\begin{eqnarray*}
  b&<&\left|\dfrac{1}{3}\left\{\int_0^r[(1+2r^{-3}s^{3})(h^{2}-1)A]ds
  \right\}\bigg|_{r=r_{0}}\right|\\[3mm]
  &\leqslant&\frac{2}{3}\left\{\int_0^r\left[
  \left(1+2r^{-3}s^{3}\right)\left|h^{2}-1\right||A|\right]ds\right\}
  \bigg|_{r=r_{0}}\\[3mm]
  &\leqslant&\frac{2}{3}A_{0}\left[\int_0^r
  \left(1+2r^{-3}s^{3}\right)ds\right]\bigg|_{r=r_{0}}=A_{0}r_{0},
\end{eqnarray*}
then $\mathcal{S}_{1}^{B}$ is nonempty.

On the other hand, $\mathcal{S}_{2}^{B}$ contains big $b$.~Here we introduce a transformation $t=br$ in $\eqref{031}$ to consider
\begin{eqnarray}
\label{035}
B(t)=t+\frac{1}{3}\int_0^t \frac{\tau^{2}}{b^{2}}\left(t\tau^{-2}-t^{-2}\tau\right)\left [2\left(h^{2}-1\right)B+\dfrac{1}{4}g'^{2}\rho^{2}b^{-2}\tau^{2}(B-A)\right]d\tau,
\end{eqnarray}
\begin{eqnarray}
\label{036}
B'(t)=1+\frac{1}{3}\int_0^t \frac{\tau^{2}}{b^{2}}\left(\tau^{-2}+t^{-3}\tau\right)\left [2\left(h^{2}-1\right)B+\dfrac{1}{4}g'^{2}\rho^{2}b^{-2}\tau^{2}(B-A)\right]d\tau.
\end{eqnarray}
In view of
\begin{eqnarray*}
\left|\dfrac{h-1}{r^{2}}\right|=
\left|\dfrac{\Phi}{r^{2}}\right|\leqslant N\left(r^{2\alpha}{R^{*}}^{2}
+{R^{*}}^{\frac{2}{1+\alpha}}\right),\,\forall r\leqslant1,
\end{eqnarray*}
we obtain that
\begin{eqnarray*}
\left|(h^{2}-1)(\tau)\right|\leqslant N\left({R^{*}}^{\frac{2}{1+\alpha}}+b^{-2\alpha}\tau^{2\alpha}{R^{*}}^{2}\right)
b^{-2}\tau^{2},\,\,\forall\tau\leqslant t\leqslant b.
\end{eqnarray*}
By the continuous of $A(r),B(r),\rho(r)$ for all $r\in[0,R]$, there exist $M_{1},M_{2},M_{3}>0$ such that $|B(t)|\leqslant M_{1},\,\,|A(t)|\leqslant M_{2},\,\,|\rho(t)|\leqslant M_{3}$ for all $t\leqslant b$, then we have
\begin{eqnarray*}
\dfrac{1}{b^{2}}\left[2(h^{2}-1)B
+\dfrac{1}{4}g'^{2}\rho^{2}b^{-2}\tau^{2}(B-A)\right]
\leqslant\dfrac{1}{b^{2}}\left[2N({R^{*}}^{\frac{2}{1+\alpha}}+{R^{*}}^{2})M_{1}
+\dfrac{1}{4}g'^{2}M_{3}^{2}(M_{1}+M_{2})\right].
\end{eqnarray*}

If $b\rightarrow\infty$, according to the integral value theorem we get
\begin{eqnarray*}
\frac{1}{3}\int_0^t\frac{\tau^{2}}{b^{2}}(t\tau^{-2}-t^{-2}\tau)\bigg[2(h^{2}-1)B
+\dfrac{1}{4}g'^{2}\rho^{2}b^{-2}\tau^{2}(B-A)\bigg]d\tau\rightarrow0.
\end{eqnarray*}

Substituting it into $\eqref{035}$, there exists a constant $R$ such
that $B$ crosses $A_{0}$ at $t_{0}=A_{0}+1$ for all $t\in[0,R]$ when
$({R^{*}}^{\frac{2}{1+\alpha}}+b^{-2\alpha}{R^{*}}^{2})b^{-2}$ is
sufficiently small.~In addition, when $b>0$ is sufficiently large,
for any $t\in(0,t_{0}]$, we have $B'(t)>0$.~ In particular,
$B'(0)=b>0$ is sufficiently large at $r=0$.~The continuity can
ensure that two sets are open. \hfill$\Box$\vskip7pt

Since the connected set $b>0$ cannot consist
of two open disjoint nonempty sets,
 there must be some value of $b$ in
  neither $\mathcal{S}_{1}^{B}$ nor $\mathcal{S}_{2}^{B}$.~
  For this value of $b$, say $b_{0}$, we have a solution
  with $0\leqslant A(r)\leqslant B(r;b_{0})\leqslant A_{0},\forall r>0$.

\begin{lemma}\label{le22}
The solution corresponding to the parameter $b_{0}$ in $\mathcal{S}^{B}_{3}$ satisfies $\lim\limits_{r\rightarrow\infty}B(r)=A_{0}$.
\end{lemma}
{\bf Proof } Since $0\leqslant A(r)\leqslant B(r)\leqslant A_{0}$ and $(rB(r))''\geqslant0$, $\forall r>0$, then $\lim\limits_{r\rightarrow\infty}
(rB)'\triangleq L\leqslant+\infty$.~If $L={+\infty}$, for convenience, we assume $G=4A_{0}$.~When $r\geqslant r_{0}$, integrating over $(r_{0},r)$ for $(rB)'>G$, we see easily that
\begin{eqnarray*}
  B(r)>\frac{r_{0}}{r}B(r_{0})+(1-\frac{r_{0}}{r})G>2A_{0},
\end{eqnarray*}
which contradicts $\lim\limits_{r\rightarrow\infty}B(r)=A_{0}$.~Using the L'Hopital's rule, we have
\begin{eqnarray*}
\lim\limits_{r\rightarrow\infty}B(r)=\lim\limits_{r\rightarrow\infty}\frac{rB(r)}{r}
=\lim\limits_{r\rightarrow\infty}(rB(r))'=L.
\end{eqnarray*}
On the one hand, $L\leqslant A_{0}$ because $B\leqslant A_{0}$.~On the other hand, applying $B\geqslant A$ and $\lim\limits_{r\rightarrow\infty}A(r)=A_{0}$, we arrive at $L\geqslant A_{0}$.~Consequently, $\lim\limits_{r\rightarrow\infty}B(r)=L=A_{0}$.
\hfill$\Box$\vskip7pt

Now we have obtained that the solution corresponding to the parameter $b_{0}$ is the solution of the boundary value problem consisting of $\eqref{004}$, $\eqref{007}$ and $\eqref{008}$.~The following will prove that the solution is unique.~To this end, we have the following lemma.

\begin{lemma}\label{le23}
The solution for the given parameter $b_{0}$ is unique.
\end{lemma}
{\bf Proof } Assume that there are two solutions
$B_{1}(r),B_{2}(r)$, and set $\Psi(r)=B_{2}(r)-B_{1}(r)$.~Then it
satisfies the boundary condition $\Psi(0)=\Psi(\infty)=0$ and the
equation
\begin{eqnarray}
\label{039}
  (r\Psi(r))''=\frac{2}{r}h^{2}(r)\Psi(r)+\frac{1}{4}g'^{2}\rho^{2}(r)r\Psi(r),
  \,\,0<r<+\infty.
\end{eqnarray}

When $r>0$, we assume
$\Psi(r)>0\,(Similarly,\,it\,\,can\,\,be\,\,seen\,\,that\,\,\Psi(r)<0)$,
therefore $\Psi'(r)=0$, $\Psi''(r)\leqslant0$ at $r=r_{0}$.~Applying
the maximum principle to $\eqref{039}$, we conclude that
$\Psi(r)\equiv0$, which contradicts the assumption that
$B_{1}(r)\neq B_{2}(r)$. \hfill$\Box$\vskip7pt

\begin{lemma}\label{le24}
If $b_{0}\in\mathcal{S}_{3}^{B}$, for $r\leqslant1$, we can find a suitably large constant $R^{*}$ such that $|r^{-1}B(r)|\leqslant R^{*}$.
\end{lemma}
{\bf Proof } According to the equation $\eqref{031}$, since $|h^{2}-1|\leqslant N\left(r^{2\alpha}{R^{*}}^{2}+{R^{*}}^{\frac{2}{1+\alpha}}\right)r^{2}$, $B-A\leqslant A_{0}$, $\rho\leqslant\rho_{0}$, $B\leqslant A_{0}$ for $r\leqslant1$, we have
\begin{eqnarray*}
  |B-b_{0}r|&\leqslant&\frac{1}{3}\int_0^r r\left[\left|2(h^{2}-1)B\right|
  +\left|\frac{1}{4}g'^{2}\rho^{2}s^{2}(B-A)\right|\right]ds\\[2mm]
  &\leqslant&\frac{1}{3}\int_0^r\left[2r^{3}N\left({R^{*}}^{2}r^{2\alpha}
  +{R^{*}}^{\frac{2}{1+\alpha}}\right)A_{0}
  +\frac{r}{4}s^{2}g'^{2}\rho_{0}^{2}A_{0}\right]ds\\[2mm]
  &\leqslant&N_{1}\left({R^{*}}^{2}r^{2\alpha+4}
  +{R^{*}}^{\frac{2}{1+\alpha}}r^{4}+r^{4}\right),
\end{eqnarray*}
where $N_{1}=\max\{\frac{2}{3}NA_{0},
\frac{1}{36}g'^{2}\rho_{0}^{2}A_{0}\}$.\\

Suppose $r=2b_{0}^{-1}A_{0}$, there is a $N_{2}=\max\{2^{2\alpha+4}A_{0}^{2\alpha+4}N_{1},2^{4}A_{0}^{4}N_{1}\}$ such that
\begin{eqnarray*}
  |B(r)-2A_{0}|
  \leqslant N_{2}\left[{R^{*}}^{2}b_{0}^{-2\alpha+4}
  +\left({R^{*}}^{\frac{2}{1+\alpha}}+1\right)b_{0}^{-4}\right].
\end{eqnarray*}

Since the left-hand side is bounded, we conclude that there exist a $C>0$ such that
\begin{eqnarray*}
N_{2}[{R^{*}}^{2}b_{0}^{-2\alpha+4}+({R^{*}}^{\frac{2}{1+\alpha}}+1)b_{0}^{-4}]
\geqslant C.
\end{eqnarray*}

That is, we have $b_{0}\leqslant
N_{3}{R^{*}}^{\frac{1}{2+\alpha}}\triangleq R^{*}$, where
$N_{3}=\max\{(C^{-1}N_{2})^{\frac{1}{2\alpha+4}},(C^{-1}N_{2})^{\frac{1}{4}}\}$.~In
view of $B(r)=b_{0}r+O(r^{2+2\alpha})(r\rightarrow0^{+})$ for all
$r\leqslant1$, we obtain that
\begin{eqnarray}
\label{038}
|r^{-1}B(r)|\leqslant|b_{0}|+|O(r^{1+2\alpha})|\leqslant(R^{*}+1)\triangleq R^{*}.
\end{eqnarray}

Thus, for $r\leq1$, we can find a suitably large constant $R^{*}$
such that $|r^{-1}B(r)|\leqslant R^{*}$.
The proof of Lemma \ref{le2} is complete.\hfill$\Box$\vskip7pt

\subsection{Proof of Lemma \ref{le3} (Existence and uniqueness of\ $\sigma(r)$)}

In order to study the solution of the boundary value problem
 be related to the equation \eqref{006} subject to the boundary
  conditions $\sigma(0)=0$ and $\sigma(\infty)=\sigma_{0}$, we
   can rewrite \eqref{006} as
\begin{eqnarray}
\label{042}
(r\sigma)''-\frac{2}{r^{2}}(r\sigma)
=\kappa\left(\sigma^{2}-\frac{m^{2}}{\kappa}\right)
(r\sigma)+\frac{2}{r^{2}}\left(h^{2}-1\right)(r\sigma).
\end{eqnarray}

Let $H=r\sigma$, we can convert \eqref{042} into the integral
equation
\begin{eqnarray}
\label{043}
H(r)=Dr^{2}+\frac{1}{3}\int_0^r \left(s^{-1}r^{2}-r^{-1}s^{2}\right)H(s)\left[\kappa\left(\frac{H^{2}}{s^{2}}
-\frac{m^{2}}{\kappa}\right)+\frac{2}{s^{2}}\left(h^{2}-1\right)\right]ds,
\end{eqnarray}
where $D$ is an arbitrary constant.

Similar to the idea in Section 3, it follows that there exists a locally continuous solution $\rho(r)=Dr+O(r^{2})\,(r\rightarrow0^{+})$ of the initial value problem consisting of $\eqref{006}$ and $\eqref{007}$.~We are interested in $D>0$, and now we define three sets as follows
\begin{eqnarray*}
\mathcal{S}^{\sigma}_{1}&=&\left\{D>0:\,\exists r_{0}>0\,\,\mbox{such that}\,\,(r\sigma(r))'|_{r=r_{0}}<0\,\,\mbox{before}\,\,H(r;D)\,\,\mbox{becomes infinite}\right\},\\
\mathcal{S}^{\sigma}_{2}&=&\left\{D>0:\,r\sigma(r)\,\,\mbox{becomes infinite before}\,\,(r\sigma(r))'\,\,\mbox{becomes zero}\right\},\\
\mathcal{S}^{\sigma}_{3}&=&\left\{D>0:\,\forall r>0,\,\sigma(r;D)\,\,\mbox{is finite},\,r\sigma(r;D)\geqslant0\right\}.
\end{eqnarray*}

It is easy to see that
$\mathcal{S}^{\sigma}_{1}\cup\mathcal{S}^{\sigma}_{2}
\cup\mathcal{S}^{\sigma}_{3}=\mathcal{S}^{\sigma},\,\mathcal{S}^{\sigma}_{1}\cap\mathcal{S}^{\sigma}_{2}
=\mathcal{S}^{\sigma}_{2}\cap\mathcal{S}^{\sigma}_{3}
=\mathcal{S}^{\sigma}_{3}\cap\mathcal{S}^{\sigma}_{1}=\emptyset$.

\begin{lemma}\label{le31}
The set $\mathcal{S}_{1}^{\sigma},\mathcal{S}_{2}^{\sigma}$ are both open and nonempty.
\end{lemma}
{\bf Proof } When $D>0$ is sufficiently small, we could assume $\sigma^{2}-\frac{m^{2}}{\kappa}=-\frac{m^{2}}{2\kappa}$ for any bounded range of $r$, then
\begin{eqnarray}
\label{044}
(r\sigma(r))''=\left(\frac{2}{r^{2}}h^{2}-\frac{m^{2}}{2}\right)(r\sigma(r)).
\end{eqnarray}

It is clear that $(r\sigma(r))'>0,r\sigma(r)>0$ for all small $r$.
~If $r>\frac{2}{m}$, $\eqref{044}$ is an oscillatory equation in the
above bounded range of $r$.~Hence, for $r>\frac{2}{m}$, there exists
a bounded range of $r$ such that $(r\sigma)'<0$.~Meanwhile,
$r\sigma$ is finite, which means $\sigma$ is finite.~Then there
exists a $r=r_{0}>0$ such that $(r\sigma)'|_{r=r_{0}}<0$ before
$H(r;D)$ becomes infinite.~Hence $\mathcal{S}_{1}^{\sigma}$ is
nonempty.

In order to proof $\mathcal{S}_{2}^{\sigma}$ is nonempty,
 we use the variable $t$ to replace $r:r=D^{-\frac{1}{2}}t$.~
 Thus   \eqref{043} becomes
\begin{eqnarray}
\label{045}
H(t)=t^{2}+\frac{1}{3}\int_0^t\left(\tau^{-1}t^{2}-\tau^{2}t^{-1}\right)
H(\tau)\left[\kappa\frac{H^{2}}{\tau^{2}}-\frac{m^{2}}{D}-
\frac{2}{\tau^{2}}\left(h^{2}-1\right)\right]d\tau.
\end{eqnarray}

Then we get
\begin{eqnarray}
\label{046}
&&H'(t)=2t+\frac{1}{3}\int_0^t\left(2\tau^{-1}t+\tau^{2}t^{-2}\right)
H(\tau)\left[\kappa\frac{H^{2}}{\tau^{2}}-\frac{m^{2}}{D}-
\frac{2}{\tau^{2}}\left(h^{2}-1\right)\right]d\tau,\\[2mm]
\label{047}
&&H''(t)=2+\frac{1}{3}\int_0^t2\tau^{-1}\left(1-\tau^{3}t^{-3}\right)
H(\tau)\left[\kappa\frac{H^{2}}{\tau^{2}}-\frac{m^{2}}{D}-
\frac{2}{\tau^{2}}\left(h^{2}-1\right)\right]d\tau\notag\\[1mm]
&&\,\,\,\,\,\,\,\,\,\,\,\,\,\,\,\,\,\,\,\,\,\,+H(\tau)\left[\kappa\frac{H^{2}}{\tau^{2}}
-\frac{m^{2}}{D}-\frac{2}{\tau^{2}}\left(h^{2}-1\right)\right].
\end{eqnarray}

According to
\begin{eqnarray*}
\left|\dfrac{\Phi}{r^{2}}\right|\leqslant N\left(r^{2\alpha}{R^{*}}^{2}
+{R^{*}}^{\frac{2}{1+\alpha}}\right),\,\forall r\leqslant1,
\end{eqnarray*}
we arrive at
\begin{eqnarray*}
\left|\frac{(h^{2}-1)}{\tau^{2}}\right|\leqslant N\left({R^{*}}^{\frac{2}{1+\alpha}}+D^{-\alpha}\tau^{2\alpha}{R^{*}}^{2}\right)D^{-1},
\,\,\forall\tau\leqslant t\leqslant D^{\frac{1}{2}}.
\end{eqnarray*}

If $D\rightarrow\infty$ in above inequation, then
\begin{eqnarray*}
\frac{m^{2}}{D}+\frac{2}{\tau^{2}}\left(h^{2}-1\right)
\leqslant\frac{m^{2}+N\left({R^{*}}^{\frac{2}{1+\alpha}}
+{R^{*}}^{2}\right)}{D}\rightarrow0.
\end{eqnarray*}

By using  \eqref{045}-\eqref{047}  and the above estimate, we obtain
that $H(t)>0,\,H'(t)>0,\,H''(t)>0$. Therefore
$\mathcal{S}_{2}^{\sigma}$ is nonempty.~Continuity can ensure that
two sets are open. \hfill$\Box$\vskip7pt

Since the connected set $D>0$ cannot consist of two open disjoint non-empty sets, there must be some value of $D$ in neither $\mathcal{S}_{1}^{\sigma}$ nor $\mathcal{S}_{2}^{\sigma}$.~For this value of $D$, say $D_{0}$, we have a solution with $\sigma(r;D)$ is finite and $r\sigma(r;D)\geqslant0,\,\forall r>0$.

\begin{lemma}\label{le32}
The solution corresponding to the parameter $D_{0}$ in $\mathcal{S}^{\sigma}_{3}$ satisfies $\lim\limits_{r\rightarrow\infty}\sigma(r)=\sigma_{0}$.
\end{lemma}
{\bf Proof } The fact that $\sigma(r;D_{0})$ is bounded at $[0,+\infty)$ follows immediately from the fact $D_{0}\in\mathcal{S}^{\sigma}_{3}$.~To see that
\begin{eqnarray*}
\lim\limits_{r\rightarrow\infty}\sigma(r)
=m\sqrt{\frac{1}{\kappa}}\triangleq\sigma_{0},
\end{eqnarray*}
we have the following two steps.

On the one hand, we have $\sigma\leqslant\sigma_{0}$.~If not, we
 assume $\sigma>\sigma_{0}$.~Next we
  claim $\lim\limits_{r\rightarrow\infty}\sigma'(r)=0$. We
   could assume $\lim\limits_{r\rightarrow\infty}\sigma'(r)=\beta>0$,
   then $\sigma'>\frac{\beta}{2}>0$ when $r$ is sufficiently
   small.~And then $\sigma(r)>\frac{\beta}{2}r+C$, which
   contradicts the finiteness of $\sigma(r)$.~Hence
    $\lim\limits_{r\rightarrow\infty}\sigma'(r)=0$.~
    Applying \eqref{042}, we observe that $\sigma''>0$
    when $r$ is sufficiently large, which contradicts $\lim\limits_{r\rightarrow\infty}\sigma'(r)=0$.

On the other hand, we have $\sigma\geqslant\sigma_{0}-\varepsilon$,
$\forall\varepsilon>0$.~If not, we assume
$\sigma<\sigma_{0}-\varepsilon$.~Similarly,
$\lim\limits_{r\rightarrow\infty}\sigma'(r)=0$ is valid.~According
to \eqref{042}, we obtain that $(r\sigma)''<0$ when $r>0$ is
sufficiently large, that is, there exist a positive constant $C$
such that $H''<-CH$.~So $D_{0}\in\mathcal{S}^{\sigma}_{1}$, which
contradicts $D_{0}\in\mathcal{S}^{\sigma}_{3}$.

With the above analysis, we must have $\sigma_{0}-\varepsilon\leqslant\sigma\leqslant\sigma_{0}$ as $r\rightarrow{+\infty}$, in other words, $\lim\limits_{r\rightarrow\infty}\sigma(r)=\sigma_{0}$.
\hfill$\Box$\vskip7pt

Finally, we will prove the uniqueness and some properties of the solution $\sigma$.~To this end, we have the following lemmas.

\begin{lemma}\label{le33}
The solution for the given parameter $D_{0}$ is unique.
\end{lemma}
{\bf Proof } Suppose otherwise that there are two solutions $\sigma_{1},\sigma_{2}$, and set $\Psi(r)=\sigma_{2}(r)-\sigma_{1}(r)$. Then the function $\Psi(r)$ satisfies the boundary condition $\Psi(0)=\Psi(\infty)=0$ and the equation
\begin{eqnarray}
\label{052}
(r\Psi(r))''&=&\left[\kappa\left(\sigma_{2}^{2}(r)-\frac{m^{2}}{\kappa}\right)
  +\frac{2}{r^{2}}h^{2}(r)
  +\kappa(\sigma_{1}^{2}(r)+\sigma_{1}(r)\sigma_{2}(r))\right](r\Psi(r))\notag\\[1mm]
&\triangleq& Q_{2}(r)(r\Psi(r)),\,\,\,0<r<+\infty.
\end{eqnarray}

Without loss of generality, we assume $\Psi(r)>0$ when $r>0$ is
sufficiently small.~In view of $\eqref{042}$, we obtain that
\begin{eqnarray}
\label{051}
  (r\sigma_{i}(r))''&=&\left[\kappa\left(\sigma_{i}^{2}(r)-\frac{m^{2}}{\kappa}\right)
  +\frac{2}{r^{2}}h^{2}(r)\right](r\sigma_{i}(r))\notag\\[1mm]
  &\triangleq& q_{i}(r)(r\sigma_{i}(r)),\,\,0<r<+\infty,\,\,(i=1,2).
\end{eqnarray}

Since $Q_{2}(r)-q_{1}(r)>0$, Applying the Sturm-Picone comparison theorem to $\eqref{052}-\eqref{051}$, we conclude that $r\sigma_{2}(r)$ have more zero points than $r\Psi(r)$ for all $r\in(0,+\infty)$.~Note that $r\sigma_{2}(r)\neq0$ for all $r\in(0,+\infty)$, then we have $\Psi(r)\neq0$ at a finite internal of $r$.~Multiplying the equation $\eqref{052}$ by $r\sigma_{2}(r)$, and the equation$\eqref{051}$ $($take $i=2)$ by $r\Psi(r)$, and then subtracting, we get \begin{equation*}
(r\Psi(r))''(r\sigma_{2}(r))-(r\sigma_{2}(r))''(r\Psi(r))
=\kappa(\sigma_{1}^{2}(r)+\sigma_{1}(r)\sigma_{2}(r))(r\Psi(r))(r\sigma_{2}(r))>0.
\end{equation*}

Then $[(r\Psi(r))'(r\sigma_{2}(r))-(r\sigma_{2}(r))'(r\Psi(r))]'>0$, that is, $(r\Psi(r))'(r\sigma_{2}(r))-(r\sigma_{2}(r))'(r\Psi(r))$ is monotonically increasing.~According to
\begin{eqnarray*}
\left[(r\Psi(r))'(r\sigma_{2}(r))-(r\sigma_{2}(r))'(r\Psi(r))\right]\bigg|_{r=0}
=\left[r^{2}(\Psi'(r)\sigma_{2}(r)-\Psi(r)\sigma_{2}'(r))\right]\bigg|_{r=0}=0,
\end{eqnarray*}
we observe that $(r\Psi(r))'(r\sigma_{2}(r))-(r\sigma_{2}(r))'(r\Psi(r))>0$.

In view of
\begin{equation*}
((\sigma_{2}^{-1}(r))\Psi(r))'
=(\sigma_{2}^{-2}(r))(\Psi'(r)\sigma_{2}(r)-\sigma_{2}'(r)\Psi(r))>0
\end{equation*}
and $(\sigma_{2}^{-1}(r))\Psi(r)>0$ at $r=0+\varepsilon$, we easily obtain $(\sigma_{2}^{-1}(r))\Psi(r)>0$ as $r\rightarrow\infty$, which contradicts \begin{eqnarray*}
  \lim\limits_{r\rightarrow\infty}\left(\frac{\Psi(r)}{\sigma_{2}(r)}\right)
  =\lim\limits_{r\rightarrow\infty}
  \left(\frac{\sigma_{2}(r)-\sigma_{1}(r)}{\sigma_{2}(r)}\right)
  =\lim\limits_{r\rightarrow\infty}\left(1-\frac{\sigma_{1}(r)}{\sigma_{2}(r)}\right)
  =0.
\end{eqnarray*}
Thus the Lemma \ref{le33} follows.
\hfill$\Box$\vskip7pt

\begin{lemma}\label{le34}
If $D_{0}\in\mathcal{S}_{3}^{\sigma}$, for $r\leqslant1$, we can find a suitably large constant $R^{*}$ such that $|r^{-1}\sigma(r)|\leqslant R^{*}$.~Moreover, if $\sigma(r)\leqslant\sigma_{0}$, then $r^{-2}H(r)$ is decreasing.
\end{lemma}
{\bf Proof } Since the equation \eqref{043} and
$H(r)>0,\sigma(r)\leqslant\sigma_{0},0\leqslant h(r)\leqslant1$ for
all $r>0$, we must have $(r^{-2}H(r))'\leqslant0$.~In other words,
$r^{-2}H(r)$ is decreasing.~To prove the other part, for all
$r\leqslant1$, we arrive at
\begin{eqnarray*}
|\sigma(r)-D_{0}r|
&\leqslant&\frac{1}{3}\int_0^r
\left|s\sigma_{0}\left(s^{-1}r-s^{2}r^{-2}\right)
\left[2N\left(s^{2\alpha}{R^{*}}^{2}
+{R^{*}}^{\frac{2}{1+\alpha}}\right)\right]\right|ds\\[2mm]
&\leqslant&\frac{1}{3}\int_0^r
\left[s^{-1}r\times s\sigma_{0}\times2N\left(s^{2\alpha}{R^{*}}^{2}
+{R^{*}}^{\frac{2}{1+\alpha}}\right)\right]ds\\[2mm]
&=&N_{1}\left(\frac{{R^{*}}^{2}}{2\alpha+1}r^{2\alpha+2}
+r^{2}{R^{*}}^{\frac{2}{1+\alpha}}\right) \leqslant
N_{2}\left(r^{2\alpha+2}{R^{*}}^{2}
+r^{2}{R^{*}}^{\frac{2}{1+\alpha}}\right).
\end{eqnarray*}
Here we have used
 $0\leqslant h(r)\leqslant1$, $H(r)\leqslant r\sigma_{0}$ and
\begin{eqnarray*}
|h^{2}(r)-1|\leqslant N\left(r^{2\alpha}{R^{*}}^{2}+{R^{*}}^{\frac{2}{1+\alpha}}\right)r^{2}
,\,\,\forall r\leqslant1\,\,(0<\alpha<k=\frac{\sqrt{3}-1}{2}),
\end{eqnarray*}
where $N_{1}=\frac{2Nr\sigma_{0}}{3},N_{2}=\max\{\frac{N_{1}}{2\alpha+1},N_{1}\}$.

Suppose $r=2\sigma_{0}D_{0}^{-1}$, there exists
$N'=\max\{2^{2+2\alpha}\sigma_{0}^{2+2\alpha}N_{2},4\sigma_{0}^{2}N_{2}\}$
such that
\begin{eqnarray*}
|\sigma(r)-2\sigma_{0}|
\leqslant N'\left(D_{0}^{-2-2\alpha}{R^{*}}^{2}+{R^{*}}^{\frac{2}{1+\alpha}}D_{0}^{-2}\right).
\end{eqnarray*}

Since the left-hand side of the above inequality is bounded, we conclude that there exist a $C>0$ such that
\begin{eqnarray*}
N'[{R^{*}}^{2}D_{0}^{-2\alpha-2}
+{R^{*}}^{\frac{2}{1+\alpha}}D_{0}^{-2}]\geqslant C.
\end{eqnarray*}

Thus, we obtain $D_{0}\leqslant N_{3}{R^{*}}^{\frac{1}{1+\alpha}}$,
where $N_{3}=\max\{(C^{-1}N')^{\frac{1}{2\alpha+2}},
(C^{-1}N')^{\frac{1}{2}}\}$.~Notice that the only positive term in
the integrand of $\eqref{043}$ is $\frac{H^{2}}{s^{2}}$ and
$\frac{H^{2}}{s^{2}}\leqslant\sigma_{0}^{2}$ is bounded.~Then for
all $\forall r\leqslant1$, we get
\begin{eqnarray}
\label{050}
\sigma(r)
&\leqslant&D_{0}r+\frac{1}{3}\int_0^r\bigg[s^{-1}r(1-s^{3}r^{-3})\times s\sigma_{0}\times\kappa\frac{H^{2}}{s^{2}}\bigg]ds
\leqslant D_{0}r+\frac{1}{3}\int_0^r\left[s^{-1}r\times s\sigma_{0}\times\kappa\sigma_{0}^{2}\right]ds\notag\\[2mm]
&=&D_{0}r+\frac{\kappa\sigma_{0}^{3}r^{2}}{3}
\leqslant N_{3}{R^{*}}^{\frac{1}{1+\alpha}}r
+\frac{\kappa\sigma_{0}^{3}}{3}r^{2}
\leqslant N\left(R^{*}r+r^{2}\right)\leqslant2NR^{*}r\triangleq R^{*}r,
\end{eqnarray}
where $R^{*}$ is a suitably large constant, $N=\max\{N_{3},\,\frac{1}{3}\kappa\sigma_{0}^{3}\}$, $0<\frac{1}{1+\alpha}<1$.
This completes the proof of the Lemma
\ref{le3}.\hfill$\Box$\vskip7pt

\subsection{Proof of Lemma \ref{le4} (Existence and uniqueness of\ $\tilde{h}(r)$)}

The proof of this subsection is similar to that of  subsection 2.2,
and only the conclusion and part of the proof are given
here.~Setting $\Phi=\tilde{h}-1$, then we may rewrite equation
$\eqref{0007}$ as
\begin{eqnarray}
\label{054}
\Phi''=\frac{2\Phi}{r^{2}}+\left(g'^{2}\sigma^{2}-B^{2}\right)\left(\Phi+1\right)
+\frac{\Phi^{3}+3\Phi^{2}}{r^{2}}.
\end{eqnarray}

From \eqref{054} and $\Phi(0)=0$ we can exhibit \eqref{054}
alternatively as the integral equation
\begin{eqnarray}
\label{055}
\Phi(r)=Er^{2}+\frac{1}{3}\int_0^r\left(r^{2}s^{-1}-r^{-1}s^{2}\right)
\left\{\left(g'^{2}\sigma^{2}-B^{2}\right)\left(\Phi+1\right)
+\frac{3\Phi^{2}+\Phi^{3}}{s^{2}}\right\}ds,
\end{eqnarray}
where $E$ is an arbitrary constant.

Then it follows that there exists a locally continuous solution
$\Phi(r)=Er^{2}+O(r^{2+2\alpha})\,(r\rightarrow0^{+})$ of the
initial value problem.~We are interested in $E<0$, and now we define
three sets as follows
\begin{eqnarray*}
\mathcal{S}_{1}^{\tilde{h}}&=&\left\{E<0:\,\Phi'(r;E)\,\,\mbox{becomes positive before}\,\,\Phi(r;E)\,\,\mbox{reaches}\,-1\right\},\\[1mm]
\mathcal{S}_{2}^{\tilde{h}}&=&\left\{E<0:\,\Phi(r;E)\,\,\mbox{crosses}\,-1\,
\,\mbox{before}\,\,\Phi'(r;E)\,\,\mbox{becomes}\,\,0\right\},\\[1mm]
\mathcal{S}_{3}^{\tilde{h}}&=&\left\{E<0:\,\forall r>0,\Phi'(r;E)\leqslant0,-1<\Phi(r;E)<0\right\}.
\end{eqnarray*}
It is clear that $\mathcal{S}^{\tilde{h}}_{1}\cup\mathcal{S}^{\tilde{h}}_{2}
\cup\mathcal{S}^{\tilde{h}}_{3}=\mathcal{S}^{\tilde{h}},
\,\mathcal{S}^{\tilde{h}}_{1}\cap\mathcal{S}^{\tilde{h}}_{2}
=\mathcal{S}^{\tilde{h}}_{2}\cap\mathcal{S}^{\tilde{h}}_{3}
=\mathcal{S}^{\tilde{h}}_{3}\cap\mathcal{S}^{\tilde{h}}_{1}=\emptyset$.

The following two lemmas are used to prove the
 existence of the solution of two-point boundary value
  problem about $\eqref{054}$, and the procedure is
  similar to Lemma \ref{le11}-Lemma \ref{le12}.
  Thus the proofs are omitted here.

\begin{lemma}\label{le41}
The set $\mathcal{S}_{1}^{\tilde{h}},\mathcal{S}_{2}^{\tilde{h}}$ are both open and nonempty.
\end{lemma}

Since the connected set $E<0$ cannot consist of two open disjoint non-empty sets, there must be some value of $E$ in neither $\mathcal{S}_{1}^{\tilde{h}}$ nor $\mathcal{S}_{2}^{\tilde{h}}$.~For this value of $E$, say $E_{0}$, we have a solution with $\Phi'(r;E)\leqslant0,\,-1<\Phi(r;E)<0$.

\begin{lemma}\label{le42}
The solution corresponding to the parameter $E_{0}$ in $\mathcal{S}_{3}^{\tilde{h}}$ satisfies $\lim\limits_{r\rightarrow\infty}\tilde{h}(r)=0$.
\end{lemma}

\begin{lemma}\label{le43}
The solution for the given parameter $E_{0}$ is unique.
\end{lemma}

The conclusion can be obtained by the Sturm-Picone comparison
theorem, and the proofs are also omitted here.

\begin{lemma}\label{le44}
If $E_{0}\in\mathcal{S}_{3}^{\tilde{h}}$, for $r\leqslant\leqslant1$, we can find $N=N(\sigma_{0},\kappa,\mu,B_{0},m,g')>0$ such that
\begin{eqnarray*}
|r^{-2}\Phi(r)|\leqslant N\left({R^{*}}^{\frac{2}{1+\alpha}}+r^{2\alpha}{R^{*}}^{2}\right),
\end{eqnarray*}
where $\Phi(r)=\tilde{h}(r)-1$.~Furthermore, we have ${\tilde{h}}'(0)=0$.
\end{lemma}
{\bf Proof } According to $\eqref{055}$, $B(s)\leqslant s^{\alpha}R^{*},\,\,\sigma(s)\leqslant s^{\alpha}R^{*},\,\,\forall s\leqslant1$, and $\left|E_{0}\right|^{1+\alpha}\leqslant N{R^{*}}^{2}$ $($because $E_{0}\in\mathcal{S}_{3}^{\tilde{h}})$, we arrive at
\begin{eqnarray}
\label{063}
\left|\frac{\Phi}{r^{2}}\right|
&\leqslant&\left|E_{0}\right|+\frac{1}{3}\int_0^r\left|s^{-1}
\left(g'^{2}+1\right)s^{2\alpha}{R^{*}}^{2}\right|ds\notag\\[1mm]
&\leqslant&\left(N{R^{*}}^{2}\right)^{\frac{1}{1+\alpha}}
+\frac{\left({g'}^{2}+1\right)}{3}\int_0^r s^{2\alpha-1}{R^{*}}^{2}ds\notag\\[2mm]
&\leqslant&N'\left({R^{*}}^{\frac{2}{1+\alpha}}+r^{2\alpha}{R^{*}}^{2}\right)
\leqslant N'\left({R^{*}}^{\frac{2}{1+\alpha}}+{R^{*}}^{2}\right)\triangleq{R_{1}^{*}},
\,\,\forall r\leqslant1,
\end{eqnarray}
where
$N'=\max\{N^\frac{1}{1+\alpha},\frac{\left(g'^{2}+1\right)}{6\alpha}\}$.~Using
the above result, we have $\Phi(r)=O(r^{2})(r\rightarrow0^{+})$. By
the defination of the derivative of the function $\Phi(r)$ at $r=0$
we obtain that
\begin{equation*}
  \Phi'(0)=\lim\limits_{r\rightarrow0^{+}}\frac{\Phi(r)-\Phi(0)}{r}
  =\lim\limits_{r\rightarrow0^{+}}\frac{O(r^{2})}{r}=0,
\end{equation*}
In other words, $\tilde{h}'(0)=0$.
The proof of Lemma \ref{le4} is complete.\hfill$\Box$\vskip7pt

\subsection{Proof of Lemma \ref{le5} (Existence and uniqueness of\ $\tilde{\rho}(r)$)}

The proof of this  subsection is similar to that of  subsection 2.4,
and only the conclusion and part of the proof are given here.~In
order to study the solution of the boundary value problem be related
to the equation $\eqref{0009}$ subject to the boundary conditions
$\tilde{\rho}(0)=0$ and $\tilde{\rho}(\infty)=\rho_{0}$, we can
rewrite $\eqref{0009}$ as
\begin{eqnarray}
\label{067}
\left(r\tilde{\rho}\right)''-\frac{1}{2r^{2}}\left(r\tilde{\rho}\right)
=-\frac{1}{4}(A-B)^{2}\left(r\tilde{\rho}\right)+\frac{\lambda}{2}\left(r\tilde{\rho}\right)
\left(\tilde{\rho}^{2}-\rho_{0}^{2}\right)
+\frac{1}{2r^{2}}(f^{2}-1)\left(r\tilde{\rho}\right).
\end{eqnarray}
Let $Q=r\tilde{\rho}$, we can convert $\eqref{067}$ into the integral equation
\begin{eqnarray}
\label{068}
Q(r)=Fr^{k+1}+\dfrac{1}{\sqrt{3}}\int_0^r \left(s^{-k}r^{k+1}-r^{-k}s^{k+1}\right)Q(s)T(s)ds,
\end{eqnarray}
where $F$ is an arbitrary constant,
\begin{eqnarray*}
T(s)=-\frac{1}{4}(A-B)^{2}+\frac{\lambda}{2}\left(\frac{Q^{2}}{s^{2}}-\rho_{0}^{2}\right)
+\frac{1}{2s^{2}}(f^{2}-1).
\end{eqnarray*}

Then it follows that there exists a locally continuous solution
$Q(r)=Fr^{k+1}+O(r^{2+k})\,(r\rightarrow0^{+})$ of the initial value
problem  consisting of $\eqref{0009}$ and $\tilde{\rho}(0)=0$. We
are interested in $F>0$, and now we define three sets as follows
\begin{eqnarray*}
\mathcal{S}^{\tilde{\rho}}_{1}&=&\left\{F>0:\,\exists r_{0}>0\,\,\mbox{such that}\,\,(r\tilde{\rho}(r))'|_{r=r_{0}}<0\,\,\mbox{before}
\,\,\tilde{\rho}(r;F)\,\,\mbox{becomes infinite}\right\},\\
\mathcal{S}^{\tilde{\rho}}_{2}&=&\left\{F>0:\,r\tilde{\rho}(r)\,\,\mbox{becomes infinite before}\,\,(r\tilde{\rho}(r))'\,\,\mbox{becomes zero}\right\},\\
\mathcal{S}^{\tilde{\rho}}_{3}&=&\left\{F>0:\,\forall r>0,\,\tilde{\rho}(r;F)\,\,\mbox{is finite},\,r\tilde{\rho}(r;F)\geqslant0\right\}.
\end{eqnarray*}

It is obvious that $\mathcal{S}^{\tilde{\rho}}_{1}\cup\mathcal{S}^{\tilde{\rho}}_{2}
\cup\mathcal{S}^{\tilde{\rho}}_{3}=\mathcal{S}^{\tilde{\rho}},
\,\mathcal{S}^{\tilde{\rho}}_{1}\cap\mathcal{S}^{\tilde{\rho}}_{2}
=\mathcal{S}^{\tilde{\rho}}_{2}\cap\mathcal{S}^{\tilde{\rho}}_{3}
=\mathcal{S}^{\tilde{\rho}}_{3}\cap\mathcal{S}^{\tilde{\rho}}_{1}=\emptyset$.

The following two lemmas are used to prove the existence of
 the solution of two-point boundary value problem about \eqref{067},
 and the procedure is similar
 to Lemma \ref{le31}-Lemma \ref{le32}.
 Thus the proof is omitted here.

\begin{lemma}\label{le51}
The set $\mathcal{S}_{1}^{\tilde{\rho}},\mathcal{S}_{2}^{\tilde{\rho}}$ are both open and nonempty.
\end{lemma}

Since the connected set $F>0$ cannot consist of two open disjoint non-empty sets, there must be some value of $F$ in neither $\mathcal{S}_{1}^{\tilde{\rho}}$ nor $\mathcal{S}_{2}^{\tilde{\rho}}$.~For this value of $F$, say $F_{0}$, we have a solution with $\tilde{\rho}(r;F)$ is finite and $r\tilde{\rho}(r;F)\geqslant0,\,\forall r>0$.

\begin{lemma}\label{le52}
If $F_{0}\in\mathcal{S}^{\tilde{\rho}}_{3}$, then $\tilde{\rho}^{2}(r;F_{0})\leqslant\rho_{0}^{2}+\frac{A_{0}^{2}}{2\lambda}$ as $r\rightarrow\infty$.
\end{lemma}
{\bf Proof } If not, we assume that $r=r_{0}$ is the first point such that $\tilde{\rho}^{2}(r)=\rho_{0}^{2}+\frac{A_{0}^{2}}{2\lambda}$.~Since $\tilde{\rho}(0)=0$, we get $\tilde{\rho}'(r_{0})\geqslant0$.~If $\tilde{\rho}'(r_{0})=0$, in view of $-\frac{1}{4}(B-A)^{2}>-\frac{A_{0}^{2}}{4}$, then ${\tilde{\rho}}''(r_{0})>0$. That is, $r_{0}$ is the minimum point of $\tilde{\rho}$, which contradicts $r_{0}$ is the maximum point of $\tilde{\rho}$.~If $\tilde{\rho'}(r_{0})>0$, it is clear that
\begin{eqnarray*}
(r\tilde{\rho})''>\left(\rho_{0}^{2}+\frac{2\lambda}{A_{0}^{2}}\right)^{\frac{1}{2}}
\left[\frac{A_{0}^{2}}{4}-\frac{1}{4}\left(B(r)-A(r)\right)^{2}\right]r\triangleq Nr,\,(N>0),
\end{eqnarray*}
because $\tilde{\rho}(\infty)>\{\rho_{0}^{2}+\frac{A_{0}^{2}}{2\lambda}\}^{\frac{1}{2}}$.~Then $\tilde{\rho}$ becomes infinite at a finite point, in other words, $F_{0}\in\mathcal{S}^{\tilde{\rho}}_{2}$, which contradicts $F_{0}\in\mathcal{S}^{\tilde{\rho}}_{3}$.~To sum up, we have $\tilde{\rho}^{2}(r;F_{0})\leqslant\rho_{0}^{2}+\frac{A_{0}^{2}}{2\lambda}$ as $r\rightarrow\infty$.
\hfill$\Box$\vskip7pt

\begin{lemma}\label{le53}
The solution corresponding to the parameter $F_{0}$ in $\mathcal{S}^{\tilde{\rho}}_{3}$ satisfies $\lim\limits_{r\rightarrow\infty}\tilde{\rho}(\infty)=\rho_{0}$.
\end{lemma}

\begin{lemma}\label{le54}
If $F_{0}\in\mathcal{S}^{\tilde{\rho}}_{3}$, $\forall\varepsilon>0$, there exist a $R(\varepsilon)$, $\forall r>R(\varepsilon)$, such that
\begin{eqnarray}
\label{0}
  \tilde{\rho}^{2}(r)>\rho_{0}^{2}-\varepsilon,
\end{eqnarray}
where $R(\varepsilon)$is independent of the choice of $\rho,A,R^{*}$ in Lemma \ref{le1}.
\end{lemma}
{\bf Proof } When
$r^{2}\geqslant2(\varepsilon\lambda)^{-1}\triangleq
R_{1}^{2}(\varepsilon)$, we assume
$\tilde{\rho}^{2}(r)\leqslant\rho_{0}^{2} -\varepsilon$.~Inserting
this result into $\eqref{067}$, we obtain that
\begin{eqnarray}
\label{000}
  Q''(r)\leqslant-\frac{\varepsilon}{2}\lambda Q(r)+\frac{1}{4}\varepsilon\lambda Q(r)=-\frac{1}{4}\varepsilon\lambda Q(r).
\end{eqnarray}

It is easy to see that $Q(r)\neq0,\forall r\in(0,+\infty)$ because
$Q(0)=0$ and $Q'(r)>0$.~On the other hand, if $r\geqslant R_{1}$,
applying the Sturm-Picone comparison theorem to
$Q_{1}''(r)=-\frac{1}{4}\varepsilon\lambda Q_{1}(r)$ and
$\eqref{000}$, we obtain that $Q(r)$ have more zero points than
$Q_{1}(r)$.~Then in view of the above equation with zeros
$2k\pi(\lambda\varepsilon)^{-\frac{1}{2}}$ $($where $k$ is an
integer$)$, there exists at least one zero point of $Q(r)$ in the
interval $(R_{1},R(\varepsilon))$ $($where $R(\varepsilon)\triangleq
R_{1}+4\pi(\lambda\varepsilon)^{-\frac{1}{2}})$, which contradicts
$Q(r)\neq0,\forall r\in(0,+\infty)$. Consequently, there is a
$r_{0}\in(R_{1},R(\varepsilon))$ such that
$\tilde{\rho}^{2}(r)>\rho_{0}^{2}-\varepsilon$.~According to the
existence of solutions of two point boundary value problem of
ordinary differential equation $\eqref{0009}$ and $Q'(r)>0$, we can
get $\tilde{\rho}^{2}(r)>\rho_{0}^{2}-\varepsilon,\,\forall
r\in(R(\varepsilon),+\infty)$. \hfill$\Box$\vskip7pt

\begin{lemma}\label{le55}
The solution for the given parameter $F_{0}$ is unique.
\end{lemma}

Following the proof of Lemma \ref{le33}, we obtain

\begin{lemma}\label{le56}
If $F_{0}\in\mathcal{S}_{3}^{\tilde{\rho}}$, for $r\leqslant1$, then
there exists a suitably large constant $R_{2}^{*}$ such that
$|r^{-k}\tilde{\rho}(r)|\leqslant R_{2}^{*}$.~Moreover, if
$\tilde{\rho}(r)\leqslant\rho_{0}$, then $r^{-k-1}Q(r)$ is
decreasing.
\end{lemma}
{\bf Proof } Imitating the proof of Lemma \ref{le34}, it can be concluded that $r^{-k-1}Q(r)$ is decreasing. To prove the other part, for all $r\leqslant1$, we arrive at
\begin{eqnarray*}
|\tilde{\rho}(r)-F_{0}r^{k}|
&\leqslant&
N_{1}\int_0^rr^{k}s^{1-k}\left(1+s^{2\alpha}{R^{*}}^{2}
+{R^{*}}^{\frac{2}{1+\alpha}}\right)ds\\[2mm]
&\leqslant&2N_{1}\left(r^{2+2\alpha}{R^{*}}^{2}
+r^{2}{R^{*}}^{\frac{2}{1+\alpha}}\right),
\end{eqnarray*}
in which we have used $0\leqslant f(r)\leqslant1$, $Q(r)\leqslant
r\left(\rho_{0}^{2}+\frac{A_{0}^{2}}{2\lambda}\right)^{\frac{1}{2}}$
and \begin{eqnarray*}
\left|f^{2}(r)-1\right|\leqslant
N\left(r^{2\alpha}{R^{*}}^{2}+{R^{*}}^{\frac{2}{1+\alpha}}\right)r^{2},
\,\,\forall r\leqslant1\,\,(0<\alpha<k=\frac{\sqrt{3}-1}{2}).
\end{eqnarray*}
Here $N_{1}=\max\{\frac{A_{0}^{2}}{4\sqrt{3}}(\rho_{0}^{2}
+\frac{A_{0}^{2}}{2\lambda})^{\frac{1}{2}},
\frac{N}{2\sqrt{3}}(\rho_{0}^{2}
+\frac{A_{0}^{2}}{2\lambda})^{\frac{1}{2}}\}$.

Suppose $r^{k}=2(\rho_{0}^{2}+\frac{A_{0}^{2}}{2\lambda})^{\frac{1}{2}}F_{0}^{-1}$,
 there exists
\begin{equation*}
N_{2}=2N_{1}\max\{2^{\frac{2(1+\alpha)}{k}}(\rho_{0}^{2}
+\frac{A_{0}^{2}}{2\lambda})^{\frac{1+\alpha}{k}},
2^{\frac{2}{k}}(\rho_{0}^{2}
+\frac{A_{0}^{2}}{2\lambda})^{\frac{1}{k}}\}
\end{equation*}
such that
\begin{eqnarray*}
\left|\tilde{\rho}(r)-2\left(\rho_{0}^{2}
+\frac{A_{0}^{2}}{2\lambda}\right)^{\frac{1}{2}}\right|
\leqslant N_{2}\left(F_{0}^{-\frac{2(1+\alpha)}{k}}{R^{*}}^{2}
+F_{0}^{-\frac{2}{k}}{R^{*}}^{\frac{2}{1+\alpha}}\right).
\end{eqnarray*}

Since the left-hand side of the above inequality is bounded, we conclude that there exist a $C>0$ such that
\begin{eqnarray*}
N_{2}(F_{0}^{-\frac{2(1+\alpha)}{k}}{R^{*}}^{2}
+F_{0}^{-\frac{2}{k}}{R^{*}}^{\frac{2}{1+\alpha}})\geqslant C.
\end{eqnarray*}

That is, we have $F_{0}\leqslant N_{3}{R^{*}}^{\frac{k}{1+\alpha}}$, where $N_{3}=\max\{(N_{2}C^{-1})^{\frac{k}{2(1+\alpha)}},
(N_{2}C^{-1})^{\frac{k}{2}}\}$.~Notice that the only positive term in the integrand of $\eqref{067}$ is $\frac{Q^{2}}{s^{2}}$ and $\frac{Q^{2}}{s^{2}}\leqslant\rho_{0}^{2}+\frac{A_{0}^{2}}{2\lambda}$ is bounded.~Then for all $\forall r\leqslant1$, we get
\begin{eqnarray}
\label{076}
\tilde{\rho}(r)
&\leqslant&F_{0}r^{k}+\frac{1}{\sqrt{3}}\int_0^r
\left[s^{-k}r^{k}\left(1-s^{2k+1}r^{-2k-1}\right)\times s\left(\rho_{0}^{2}+\frac{A_{0}^{2}}{2\lambda}\right)^{\frac{1}{2}}
\times\frac{\lambda Q^{2}}{2s^{2}}\right]ds\notag\\[2mm]
&\leqslant&F_{0}r^{k}+\frac{1}{\sqrt{3}}\int_0^r r^{k}s^{1-k}\frac{\lambda}{2}\left(\rho_{0}^{2}
+\frac{A_{0}^{2}}{2\lambda}\right)^{\frac{3}{2}}ds
\leqslant N_{3}{R^{*}}^{\frac{k}{1+\alpha}}r^{k}
+\frac{\lambda r^{2}}{2\sqrt{3}}\left(\rho_{0}^{2}
+\frac{A_{0}^{2}}{2\lambda}\right)^{\frac{3}{2}}\notag\\[2mm]
&\leqslant&N_{4}\left({R^{*}}^{\frac{k}{1+\alpha}}r^{k}+r^{2}\right)
\leqslant2N_{4}R^{*}r^{k}\triangleq R_{2}^{*}r^{k},
\end{eqnarray}
where $R_{2}^{*}$ is a suitably large constant, $N_{4}=\max\{N_{3}\,,\frac{\lambda}{2\sqrt{3}}(\rho_{0}^{2}
+\frac{A_{0}^{2}}{2\lambda})^{\frac{3}{2}}\}$, $0<\frac{1}{1+\alpha}<1$. That is, $\forall r\leqslant1$, we arrive at $|r^{-k}\tilde{\rho}|\leqslant R_{2}^{*}$.
The proof of Lemma \ref{le5} is finished. \hfill$\Box$\vskip7pt

\subsection{Proof of Lemma \ref{le6} (Existence and uniqueness of\
$\tilde{A}(r)$)}

The proof of this subsection is similar to that of subsection 2.3,
and only the conclusion and part of the proof are given
here.~Firstly, we  rewrite \eqref{0011} as
\begin{eqnarray}
\label{080}
r^{2}\tilde{A}''+2r\tilde{A}'-2\tilde{A}
=2\left(f^{2}-1\right)\tilde{A}-\frac{1}{4}g^{2}\rho^{2}r^{2}\left(B-\tilde{A}\right).
\end{eqnarray}

Then the differential equation \eqref{080} can be transformed into
the integral equation form
\begin{eqnarray}
\label{081}
\tilde{A}(r)=ar+\frac{1}{3}\int_0^r s^{2}\left(rs^{-2}-r^{-2}s\right)\left [2\left(f^{2}-1\right)\tilde{A}-\frac{1}{4}g^{2}\rho^{2}s^{2}
\left(B-\tilde{A}\right)\right]ds,
\end{eqnarray}
for arbitrary constant $a$, which can be determined by
 Picard iteration.
  There exists a locally continuous solution to
  the initial value problem for $r$
  sufficiently small. Then we obtain solution
   $\tilde{A}(r)=ar+O(r^{2})$ which
  continuously depends on the parameter $a$.

Next, we introduce three sets
\begin{eqnarray*}
\mathcal{S}^{\tilde{A}}_{1}&=&\left\{a>0:\,(r\tilde{A}(r))'\,\,\mbox{becomes negative before}\,\,\tilde{A}(r;a)=B(r)\right\},\\[1mm]
\mathcal{S}^{\tilde{A}}_{2}&=&\left\{a>0:\,\tilde{A}(r;a)\,\,\mbox{crosses}\,\,B(r)
\,\,\mbox{before}\,\,(r\tilde{A}(r))'=0\,\right\},\\[1mm]
\mathcal{S}^{\tilde{A}}_{3}&=&\left\{a>0:\,\forall r>0,B(r;b)\leqslant A_{0},B(r;b)\geqslant A(r)\right\}.
\end{eqnarray*}
{\bf Case (i)} $(r\tilde{A}(r))'=0,$\,$\tilde{A}(r)=B(r)$ at $r=r_{0}$, while $(r\tilde{A}(r))'\geqslant0,$\,$\tilde{A}(r)\leqslant B(r)$ for all $r<r_{0}$.\\[1mm]
We claim that case (i) is not valid.~Inserting $\tilde{A}(r_{0})=B(r_{0})$ into $\eqref{0011}$, then we easily have $(r\tilde{A}(r))''|_{r=r_{0}}>0$.~In addition, it is clear that $(r\tilde{A}(r))'<0,\,\forall r\in(r_{0}-\delta,r_{0})$ because $(r\tilde{A}(r))'|_{r=r_{0}}=0$, which contradicts $(r\tilde{A}(r))'\geqslant0$ for all $r<r_{0}$.\\

Obviously,
$\mathcal{S}^{\tilde{A}}_{1}\cup\mathcal{S}^{\tilde{A}}_{2}
\cup\mathcal{S}^{\tilde{A}}_{3}=\mathcal{S}^{\tilde{A}},\,
\mathcal{S}^{\tilde{A}}_{1}\cap\mathcal{S}^{\tilde{A}}_{2}
=\mathcal{S}^{\tilde{A}}_{2}\cap\mathcal{S}^{\tilde{A}}_{3}
=\mathcal{S}^{\tilde{A}}_{3}\cap\mathcal{S}^{\tilde{A}}_{1}=\emptyset$.

\begin{lemma}\label{le61}
The set $\mathcal{S}_{1}^{\tilde{A}},\mathcal{S}_{2}^{\tilde{A}}$ are both open and nonempty.
\end{lemma}
{\bf Proof } Firstly, $\mathcal{S}_{1}^{\tilde{A}}$ contains small $a$.~Inserting $a=0$ into the equation $\eqref{081}$, we obtain the equations as follows
\begin{eqnarray}
\label{082}
\tilde{A}(r)=\frac{1}{3}\int_0^r s^{2}\left(rs^{-2}-r^{-2}s\right)\left [2\left(f^{2}-1\right)\tilde{A}-\frac{1}{4}g^{2}\rho^{2}s^{2}(B-\tilde{A})\right]ds.
\end{eqnarray}
In fact, when $r>0$ is sufficiently small, $\tilde{A}(r)<0$ and $r\tilde{A}(r)$ is decreasing because of
\begin{eqnarray*}
  (r\tilde{A}(r))'=\int_0^r\left[\frac{2}{s}f^{2}\tilde{A}
  -\frac{1}{4}g^{2}\rho^{2}s\left(B-\tilde{A}\right)\right]ds<0.
\end{eqnarray*}
When $a>0,r>0$ is sufficiently small, since $r\tilde{A}(r)=ar^{2}+O(r^{3})(r\rightarrow0^{+})$, we have $r\tilde{A}(r)>0$.~As $r$ increases, when $r>a$ is sufficiently small, the dominant term becomes $O(r^{3})$, and by the negative coefficient of $r^{3}$ we conclude that $r\tilde{A}(r)<0$.~Therefore, $(r\tilde{A}(r))'<0$ when $r$ is sufficiently small.~Meanwhile, $B(r)>\tilde{A}(r)$.~So $\mathcal{S}_{1}^{\tilde{A}}$ is nonempty.

The proof that the set $\mathcal{S}_{2}^{\tilde{A}}$ is a nonempty
 set is similar to Lemma \ref{le21} and will be omitted here.
 ~The continuity can ensure that two sets are open.
\hfill$\Box$\vskip7pt

Since the connected set $a>0$ cannot consist of two open disjoint nonempty sets, there must be some value of $a$ in neither $\mathcal{S}_{1}^{\tilde{A}}$ nor $\mathcal{S}_{2}^{\tilde{A}}$.~For this value of $a$, say $a_{0}$, we have a solution with $B(r;b)\leqslant A_{0},B(r;b)\geqslant A(r),\forall r>0$.

\begin{lemma}\label{le62}
The solution corresponding to the parameter $a_{0}$ in $\mathcal{S}^{\tilde{A}}_{3}$ satisfies $\lim\limits_{r\rightarrow\infty}\tilde{A}(r)=A_{0}$.
\end{lemma}
{\bf Proof } Since $a_{0}\in\mathcal{S}^{\tilde{A}}_{3}$, then
$(r\tilde{A}(r))'\geqslant0,\,\tilde{A}(r)\leqslant B(r)\leqslant
A_{0}$ for all $r>0$.~In view of $(r\tilde{A}(r))|_{r=0}=0$ and
$(r\tilde{A}(r))'>0$, we obtain that $r\tilde{A}\geqslant0$ for all
$r$.~Hence $0\leqslant\tilde{A}\leqslant B\leqslant A_{0}$ for all
$r$, that is, $\tilde{A}(r)$ is bounded on $[0,+\infty)$.

We claim that $B(r)\leqslant\tilde{A}(r)$ when $r>0$ is sufficiently
large.~If not, inserting $B(r)>\tilde{A}(r)$ into $\eqref{0011}$, we
obtain that
\begin{eqnarray*}
\left(r\tilde{A}(r)\right)''=\left[\frac{2}{r^{2}}f^{2}(r)\tilde{A}(r)
-\frac{1}{4}g^{2}\tilde{\rho}^{2}(r)r\left(B(r)-\tilde{A}(r)\right)\right]r<-Kr.
\end{eqnarray*}

That is, $(r\tilde{A}(r))''\thicksim-Kr,\,(r\tilde{A}(r))'\thicksim-Kr^{2}(r\rightarrow\infty)$, which contradicts $(r\tilde{A}(r))'\geqslant0$. Combining $B(r)\leqslant\tilde{A}(r)$ and $\tilde{A}(r)\leqslant B(r)$ as $r\rightarrow+\infty$, we arrive at
\begin{eqnarray*}
\lim\limits_{r\rightarrow\infty}\tilde{A}(r)
=\lim\limits_{r\rightarrow\infty}B(r)=A_{0}.
\end{eqnarray*}
\hfill$\Box$\vskip7pt

\begin{lemma}\label{le63}
The solution for the given parameter $a_{0}$ is unique.
\end{lemma}

The conclusion can be obtained by using the extremum principle, and the proof will be omitted here.

\begin{lemma}\label{le64}
If $a_{0}\in\mathcal{S}_{3}^{\tilde{A}}$, for $r\leqslant1$, we can find a suitably large constant $R_{3}^{*}$ such that $|r^{-1}\tilde{A}(r)|\leqslant R_{3}^{*}$.~Moreover, $r^{-1}\tilde{A}(r)$ is decreasing.
\end{lemma}
{\bf Proof } Imitating the proof of Lemma \ref{le34}, it can be concluded that $r^{-1}\tilde{A}(r)$ is decreasing. To prove the other part, for all $r\leqslant1$, we arrive at
\begin{eqnarray*}
  \left|\tilde{A}(r)-a_{0}r\right|
  &\leqslant&\frac{1}{3}\int_0^r |2r(f^{2}-1)\tilde{A}|
  +\left|\frac{r}{4}g^{2}\rho^{2}s^{2}\left(B-\tilde{A}\right)\right|ds\\[2mm]
  &\leqslant&\frac{1}{3}\int_0^r r\left [2Na_{0}s^{3}\left({R^{*}}^{2}s^{2\alpha}+{R^{*}}^{\frac{2}{1+\alpha}}\right)
  +\frac{1}{2}g^{2}s^{2k+2}{R^{*}}^{2}A_{0}\right]ds\\[2mm]
  &\leqslant&N_{1}\left[{R^{*}}^{2}\left(a_{0}r^{2\alpha+5}+A_{0}r^{2k+4}\right)
  +{R^{*}}^{\frac{2}{1+\alpha}}a_{0}r^{5}\right]
\end{eqnarray*}
because $|B(r)-\tilde{A}(r)|\leqslant2A_{0}$, $\rho(r)\leqslant\rho_{0}$, $\tilde{A}(r)\leqslant a_{0}r$ and
\begin{eqnarray*}
\left|f^{2}(r)-1\right|\leqslant N\left(r^{2\alpha}{R^{*}}^{2}+{R^{*}}^{\frac{2}{1+\alpha}}\right)r^{2},
\,\,\forall r\leqslant1,
\end{eqnarray*}
where $N_{1}=\max\{\frac{2}{3}N,\frac{1}{6}g^{2}\}$.

Suppose $r=2a_{0}^{-1}A_{0}$, there is a
\begin{equation*}
N_{2}=\max\{2^{2\alpha+5}A_{0}^{2\alpha+5}N_{1},
2^{2k+4}A_{0}^{2k+5}N_{1},2^{5}A_{0}^{5}N_{1}\}
\end{equation*}
such that
\begin{eqnarray*}
  |\tilde{A}(r)-2A_{0}|
  \leqslant N_{2}\left[{R^{*}}^{2}\left(a_{0}^{-2\alpha-4}+a_{0}^{-2k-4}\right)
  +{R^{*}}^{\frac{2}{1+\alpha}}a_{0}^{-5}\right].
\end{eqnarray*}

Since the left-hand side of the above inequality is bounded, we conclude that there exist a $a>0$ such that
\begin{eqnarray*}
N_{2}[{R^{*}}^{2}(a_{0}^{-2\alpha-4}+a_{0}^{-2k-4})
  +{R^{*}}^{\frac{2}{1+\alpha}}a_{0}^{-5}]\geqslant C.
\end{eqnarray*}

That is, we have $a_{0}\leqslant N_{3}{R^{*}}^{\frac{1}{2+\alpha}}\triangleq R_{3}^{*}$, where $N_{3}=\max\{(C^{-1}N_{2})^{\frac{1}{2\alpha+4}},
(C^{-1}N_{2})^{\frac{1}{2k+4}},(C^{-1}N_{2})^{\frac{1}{5}}\}$. In view of $r^{-1}\tilde{A}\leqslant a_{0}$ for all $r\leqslant1$, there exist a suitably large constant $R_{3}^{*}$ such that
\begin{eqnarray}
\label{086}
r^{-1}\tilde{A}(r)\leqslant a_{0}\leqslant R_{3}^{*}.
\end{eqnarray}

The proof of Lemma \ref{le6} is complete. \hfill$\Box$\vskip7pt

\section{Proof of Theorem \ref{th1} and \ref{th2}}

In this section, we complete the proof of Theorems
\ref{th1}-\ref{th2} by the Schauder fixed point theorem and the
extremum principle.

Firstly, we define the space $\mathscr{B}$ as follows
\begin{eqnarray*}
\mathscr{B}=\left\{\left(\rho,A,\Phi\right)|\rho,A,\Phi\in C\left([0,+\infty)\right),r^{-\alpha}\left(1+r^{\alpha}\right)\rho,
r^{-\alpha}\left(1+r^{\alpha}\right)A,r^{-\alpha}\left(1+r^{\alpha}\right)\Phi\,\,
\mbox{are bounded}\right\}
\end{eqnarray*}
with
\begin{eqnarray*}
\left\|\left(\rho,A,\Phi\right)\right\|_{\mathscr{B}}
=\sup_{r\in[0,+\infty)}\left\{|r^{-\alpha}\left(1+r^{\alpha}\right)\rho|
+|r^{-\alpha}\left(1+r^{\alpha}\right)A|
+|r^{-\alpha}\left(1+r^{\alpha}\right)\Phi|\right\}.
\end{eqnarray*}


Next, we define the mapping
$F{\rm:}\,(\rho,A,\Phi)\rightarrow(\tilde{\rho},\tilde{A},\tilde{\Phi})$
on $\mathscr{B}$.~In the following, we will demonstrate that it is
continuous.

\begin{lemma}\label{le71}
$F{\rm:}\,(\rho,A,\Phi)\rightarrow(\tilde{\rho},\tilde{A},\tilde{\Phi})$ is continuous on $\mathscr{B}$.
\end{lemma}
{\bf Proof } In order to obtain this conclusion, we shall show that if $\left(\rho_{1},A_{1},\Phi_{1}\right),\left(\rho_{2},A_{2},\Phi_{2}\right)\in \mathscr{B}$, then when $\left\|\left(\rho_{1}-\rho_{2},A_{1}-A_{2},
\Phi_{1}-\Phi_{2}\right)\right\|_{\mathscr{B}}\rightarrow0$, we have \begin{eqnarray*}
\left\|F\left(\rho_{1},A_{1},\Phi_{1}\right)
-F\left(\rho_{2},A_{2},\Phi_{2}\right)\right\|_{\mathscr{B}}\rightarrow0.
\end{eqnarray*}
Since $F\left(\rho_{1}\left(\infty\right)\right)
=F\left(\rho_{2}\left(\infty\right)\right)=\rho_{0}$, then for any $\varepsilon>0$ there is $R_{0}>0$ such that
\begin{eqnarray}
\label{090}
 \sup_{r\in[R_{0},+\infty)}\left|r^{-\alpha}\left(1+r^{\alpha}\right)
 \left[F\left(\rho_{1}\left(\infty\right)\right)
 -F\left(\rho_{2}\left(\infty\right)\right)\right]\right|<\varepsilon.
\end{eqnarray}

Likewise, we have
\begin{eqnarray}
\label{091}
&&\sup_{r\in[R_{0},+\infty)}\left|r^{-\alpha}\left(1+r^{\alpha}\right)
\left[F\left(A_{1}\left(\infty\right)\right)
-F\left(A_{2}\left(\infty\right)\right)\right]\right|<\varepsilon,\\[1mm]
&&\sup_{r\in[R_{0},+\infty)}\left|r^{-\alpha}\left(1+r^{\alpha}\right)
\left[F\left(\Phi_{1}\left(\infty\right)\right)
-F\left(\Phi_{2}\left(\infty\right)\right)\right]\right|<\varepsilon.
\end{eqnarray}

For the $\varepsilon$ given above, there is $\delta>0$ such that
\begin{equation*}
\sup_{r\in(0,\delta]}\left|r^{-\alpha}\left(1+r^{\alpha}\right)
\left(\rho_{1}-\rho_{2}\right)\right|
\rightarrow0,\,\left(n\rightarrow\infty\right)
\end{equation*}
and $F\left(\rho_{1}(0)\right)=F\left(\rho_{2}(0)\right)=0$, then
\begin{eqnarray}
\label{092}
 \sup_{r\in(0,\delta]}\left|r^{-\alpha}\left(1+r^{\alpha}\right)
 \left[F\left(\rho_{1}(0)\right)-F\left(\rho_{2}(0)\right)\right]\right|
 <\varepsilon.
\end{eqnarray}
Similarly, we get
\begin{eqnarray}
\label{093}
\sup_{r\in(0,\delta]}\left|r^{-\alpha}(1+r^{\alpha})
\left[F\left(A_{1}(0)\right)-F\left(A_{2}(0)\right)\right]\right|
 <\varepsilon,\\[1mm]
\sup_{r\in(0,\delta]}\left|r^{-\alpha}(1+r^{\alpha})
\left[F\left(\Phi_{1}(0)\right)-F\left(\Phi_{2}(0)\right)\right]\right|
 <\varepsilon.
\end{eqnarray}

In view of the continuity of $F$ at $[\delta,R_{0}]$, thus $F$ is continuous on
$\mathscr{B}$.
\hfill$\Box$\vskip7pt

After that, we define the nonempty bounded closed convex subset
$\mathcal{S}$ of $\mathscr{B}$ by
\begin{eqnarray*}
&&\mathcal{S}=\big\{(\rho,A,\Phi)\in\mathscr{B}\,\big|~{\rm(1)}\,\,\mbox{for}\,\,r
\leqslant1,\,\,\mbox{we have}\,\,|r^{-k}\rho|\leqslant R^{*},|r^{-1}A|\leqslant R^{*},|r^{-2}\Phi|\leqslant R^{*};\\[1mm]
&&~~~~~~~~~~~~~~~~~~~~~~~~~~~~
{\rm(2)}\,\,r\rho(r)\,,\,rA(r)\,\,\mbox{is increasing},
\,\,r^{-1}\Phi(r)\,\,\mbox{is decreasing};\\[1mm]
&&~~~~~~~~~~~~~~~~~~~~~~~~~~~~~~~~~
\mbox{if}\,\,r\leqslant\frac{\rho_{0}}{R^{*}},\,\,
\mbox{then}\,\,r^{-k}\rho(r)\,\,\mbox{is decreasing};\\[1mm]
&&~~~~~~~~~~~~~~~~~~~~~~~~~~~~
{\rm(3)}\,\,\rho^{2}\leqslant\rho_{0}^{2}+\frac{A_{0}^{2}}{2\lambda},\,
\frac{1}{4}g^{2}\rho^{2}\geqslant A^{2},\,g'^{2}\sigma^{2}\geqslant
B^{2},\,A\leqslant A_{0},\,|\Phi|\leqslant1,\\[1mm]
&&~~~~~~~~~~~~~~~~~~~~~~~~~~~~~~~~~
\mbox{if}\,\,r\geqslant R(\varepsilon,F_{0}),\,\,\mbox{then}\,\,\frac{1}{4}g^{2}\rho^{2}
\geqslant\frac{1}{2}(\frac{1}{4}g_{0}^{2}\rho_{0}^{2}+A_{0}^{2});\\[1mm]
&&~~~~~~~~~~~~~~~~~~~~~~~~~~~~
{\rm(4)}\,\,A(\infty)=A_{0},\,\rho(\infty)=\rho_{0},\,\Phi(\infty)=-1\big\}\,,
\end{eqnarray*}
where $R^{*}=\max\{R_{1}^{*},\,R_{2}^{*},\,R_{3}^{*}\}$, $R(\varepsilon,F_{0})$ is $R(\varepsilon)$ in Lemma \ref{le54}..~It is straightforward that the set $\mathcal{S}$ is indeed nonempty, bounded, closed and convex.

Finally, we show that the mapping $F$ is compact which maps
$\mathcal{S}$ into itself.
\begin{lemma}\label{le72}
For the mapping $F{\rm:}\,(\rho,A,\Phi)\rightarrow(\tilde{\rho},\tilde{A},\tilde{\Phi})$,\\[2mm]
{\rm(1)} The mapping $F$ takes $\mathcal{S}$ into itself;\\[2mm]
{\rm(2)} $F$ is compact.
\end{lemma}
{\bf Proof } Part ${\rm(1)}$ is already established in Lemmas \ref{le1}-\ref{le6} or their proofs.~Next we concentrate on the proof of Part ${\rm(2)}$.~Firstly, from the fact that the mapping $F$ is continuous on $\mathscr{B}$ and $\mathcal{S}$ is the subset of $\mathscr{B}$, it follows that $F$ is continuous on $\mathcal{S}$.~Secondly, to show $F$ is compact, we will demonstrate that the mapping $F$ maps an arbitrary bounded set in $\mathcal{S}$ to a column-compact set.~That is if $\left\{\rho_{n}(r),A_{n}(r),\Phi_{n}(r)\right\}$ is the arbitrary bounded sequences in $\mathcal{S}$, then we must prove that $\{\tilde{\rho_{n}}(r),\tilde{A_{n}}(r),\tilde{\Phi_{n}}(r)\}$ have convergent sub-sequence in $\mathcal{S}$.~Take $\{\tilde{A_{n}}(r)\}$ as an example.

Using Lemma \ref{le6} and $\{\tilde{A_{n}}(r)\}\in\mathcal{S}$, it is straightforward to show
\begin{eqnarray*}
\tilde{A_{n}}'(r)=a+\frac{1}{3}\int_0^r s^{2}\left(s^{-2}+2r^{-3}s\right)
\left\{2\left(f_{n}^{2}-1\right)\tilde{A_{n}}-\frac{1}{4}g^{2}{\tilde{\rho_{n}}}^{2}
\left(B_{n}-\tilde{A_{n}}\right)\right\}ds
\end{eqnarray*}
is bounded on any inner closed subinterval of $(0,+\infty)$.~In other words, there is $L_{1}>0$ such that $\|\tilde{A_{n}}'(r)\|_{\mathscr{B}}\leqslant L_{1}$ on any inner closed subinterval of $(0,+\infty)$, where $L_{1}$ have no connection with $n$.~Applying the mean value theorem, for any $\varepsilon>0,\,r_{1},r_{2}\in[\delta,R]\subset(0,+\infty)$ there exist $\delta=\frac{\varepsilon}{L_{1}+1}>0$ such that $|\tilde{A_{n}}(r_{1})-\tilde{A_{n}}(r_{2})|
=|\tilde{A_{n}}'(\xi)||r_{1}-r_{2}|<\varepsilon$ when $|r_{1}-r_{2}|<\delta$, where $\xi$ is between $r_{1}$ and $r_{2}$.~Hence the equicontinuity of $\{\tilde{A_{n}}(r)\}$ is apparently established.~It is clear that $\{\tilde{A_{n}}(r)\}$ is uniformly bounded because $\{\tilde{A_{n}}(r)\}\in\mathcal{S}$.~According to the Arzela-Ascoli theorem, there is a subsequence of $\{\tilde{A_{n}}(r)\}$~$($denoted as $\{\tilde{A_{n}}\})$ and the continuous function $\tilde{A}(r)$, which uniformly converges to $\tilde{A}(r)$ in any compact subinterval of $(0,+\infty)$ $($denoted as $[\delta,R])$, therefore
\begin{eqnarray}
\label{095}
\sup_{r\in[\delta,R]}\big|r^{-\alpha}(1+r^{\alpha})(\tilde{A_{n}}(r)-\tilde{A}(r))\big|
\leqslant\left(1+\sigma^{-\alpha}\right)\sup_{r\in[\delta,R]}
\big|\tilde{A_{n}}(r)-\tilde{A}(r)\big|
\rightarrow0\,(n\rightarrow\infty).
\end{eqnarray}

We only need to prove that when $r\rightarrow0$ and
$r\rightarrow+\infty$$($that is in interval $(0,\delta)$ and
$(R,+\infty)$$)$, such that if $n\rightarrow+\infty$
\begin{equation*}
\|\tilde{A_{n}}(r)-\tilde{A}(r)\|_{\mathscr{B}}\rightarrow0.
\end{equation*}

As $r\rightarrow0$$\,($that is in interval $(0,\delta)$$)$,~we can get $\tilde{A_{n}}(r)\leqslant rR^{*},$\,$\tilde{A}(r)\leqslant rR^{*}$ for all $r\leqslant1$ because of $\{\tilde{A_{n}}(r)\},\tilde{A}(r)\in\mathcal{S}$.~Then, for the $\varepsilon>0$ given above, we can choose a sufficiently small $\delta=\left(\frac{\varepsilon}{4R^{*}}\right)^{\frac{1}{1-\alpha}}$ such that
\begin{eqnarray*}
\sup_{r\in(0,\delta)}|r^{-\alpha}(1+r^{\alpha})(\tilde{A_{n}}(r)-\tilde{A}(r))|
<2\sup_{r\in(0,\delta)}|r^{1-\alpha}r^{-1}(\tilde{A_{n}}(r)-\tilde{A}(r))|
<4\delta^{1-\alpha}R^{*}\leqslant\varepsilon.
\end{eqnarray*}

With $\delta>0$ fixed, we can then find $n$ sufficiently large that
\begin{eqnarray}
\label{096}
\sup_{r\in(0,\delta)}\big|r^{-\alpha}(1+r^{\alpha})(\tilde{A_{n}}(r)-\tilde{A}(r))
\big|\rightarrow0\,(n\rightarrow\infty).
\end{eqnarray}

As $r\rightarrow+\infty$$\,($that is in interval $(R,+\infty)$$)$, we want to show that
\begin{equation*}
\sup_{r\in(R,\infty)}\big|r^{-\alpha}(1+r^{\alpha})(\tilde{A_{n}}(r)-\tilde{A}(r))\big|
<\varepsilon
\end{equation*}
for $\forall\varepsilon>0$.~It is only necessary to prove that $|\tilde{A_{n}}(r)-\tilde{A}(r)|<\varepsilon$.~From $\eqref{0011}$, we have
\begin{eqnarray*}
(r\tilde{A_{n}(r)})''=\dfrac{2}{r}f_{n}^{2}(r)\tilde{A_{n}}(r)
-\dfrac{1}{4}g^{2}\tilde{\rho_{n}}^{2}(r)r(B_{n}(r)-\tilde{A_{n}}(r))
\leqslant\frac{2}{r}f_{n}^{2}(r)\tilde{A_{n}}(r).
\end{eqnarray*}
Then, integrating the above inequality over $(r,\infty)$, we have
\begin{eqnarray*}
(s\tilde{A_{n}})'(\infty)-(s\tilde{A_{n}})'(r)
\leqslant\int_r^{+\infty}\frac{2}{s}f_{n}^{2}\tilde{A_{n}}ds.
\end{eqnarray*}
In view of
$\tilde{A_{n}}(\infty)=A_{0},$\,$\tilde{A_{n}}'(\infty)=O(r^{-2})$,
we obtain that
\begin{eqnarray*}
|A_{0}-\tilde{A_{n}}(r)|\leqslant r\tilde{A_{n}}'(r)
+\int_r^{+\infty}\frac{2}{s}f_{n}^{2}\tilde{A_{n}}ds.
\end{eqnarray*}
Since $f_{n}(r)=O(e^{-\kappa(1-\varepsilon)r}),\tilde{A_{n}}(r)=A_{0}+O(r^{-1}),
\tilde{A_{n}}'(r)=O(r^{-2})$, therefore we can get
$|A_{0}-\tilde{A_{n}}(r)|\rightarrow0(n\rightarrow\infty)$ as $r\rightarrow\infty$, that is, $|\tilde{A}(r)-\tilde{A_{n}}(r)|\rightarrow0\,(n\rightarrow\infty)$.~Thus \begin{eqnarray}
\label{097}
\sup_{r\in(R,+\infty)}|r^{-\alpha}(1+r^{\alpha})(\tilde{A_{n}}(r)-\tilde{A}(r))|
\rightarrow0\,(n\rightarrow\infty)\,.
\end{eqnarray}

In summary, we have proved for any $r$, if $n\rightarrow\infty$,then $\|\tilde{A_{n}}(r)-\tilde{A}(r)\|_{\mathscr{B}}\rightarrow0$.

Evidenced by the same token, it is easy to show that
$\|\tilde{\rho_{n}}(r)-\tilde{\rho}(r)\|_{\mathscr{B}}
\rightarrow0\,(n\rightarrow\infty)$,
$\|\tilde{\Phi_{n}}(r)-\tilde{\Phi}(r)\|_{\mathscr{B}}
\rightarrow0\,(n\rightarrow\infty)$ as well. \hfill$\Box$\vskip7pt

In conclusion, according to Lemma \ref{le71}-\ref{le72},
the mapping $F$ satisfies the conditions of the Schauder
 fixed point theorem.~Applying the Schauder fixed point
  theorem, we obtain that $(f(r),\rho(r),A(r),B(r),h(r),\sigma(r))$
  is the solution of the two-point boundary
  value problem \eqref{001}-\eqref{008} for
  the system of nonlinear ordinary differential
  equations.~Thus the Theorem \ref{th1} is proved.

Next we show some properties of the solutions obtained above.\\[2mm]
{\bf Proof of Theorem \ref{th2}}

Define the comparison function
$\eta(r)=Ce^{-\zeta(1-\varepsilon)r}$, where
$\zeta=\sqrt{g'^{2}\sigma_{0}^{2}-A_{0}^{2}}$, $C>0$ is a constant
to be chosen later, $\varepsilon>0$ is sufficiently small.~From the
equation $\eqref{005}$ and the property $h(r)>0$ we obtain that for
any $\varepsilon>0$, there is a sufficiently large
$r_{\varepsilon}>0$ so that
\begin{eqnarray}
\label{105}
\left(h-\eta\right)''&=&\left(g'^{2}\sigma^{2}-B^{2}\right)h
+\frac{1}{r^{2}}h\left(h^{2}-1\right)
-\zeta^{2}\left(1-\varepsilon\right)^{2}\eta\notag\\[1mm]
&=&\zeta^{2}\left(1-\varepsilon\right)^{2}\left(h-\eta\right)
+\zeta^{2}\left[1-\left(1-\varepsilon\right)^{2}\right]h
+\frac{1}{r^{2}}h\left(h^{2}-1\right)\notag\\[1mm]
&\geqslant&\zeta^{2}\left(1-\varepsilon\right)^{2}\left(h-\eta\right),\,r>r_{\varepsilon}.
\end{eqnarray}
Taking $C>0$ be large enough to make $(h-\eta)(r_{\varepsilon})\leqslant0$.~Thus, in view of this and the boundary condition $(h-\eta)(r)\rightarrow0\left(r\rightarrow\infty\right)$, we obtain by
applying the maximum principle theorem in $\eqref{105}$ the result $0<h\leqslant\eta=Ce^{-\zeta(1-\varepsilon)r},\,r>r_{\varepsilon}$ as expected in $\eqref{04}$.

As $r\rightarrow+\infty$, the asymptotic estimate of the function $f$ can be proved by imitating $h$.~In other words, the estimate for $f$ in $\eqref{01}$ is established.

Next, set $\tau(r)=r(B-A)$, then $\tau>0$.~By virtue of the equations $\eqref{003}$ and $\eqref{004}$, we have
\begin{equation*}
  \tau''=\frac{1}{4}\left(g'^{2}+g^{2}\right)\rho^{2}\tau
+\frac{2}{r}\left(h^{2}B-f^{2}A\right),
\end{equation*}
then $\tau''=\frac{1}{4}(g^{2}+g'^{2})\rho_{0}^{2}\tau$ as $r\rightarrow\infty$.~Now let $\eta(r)=Ce^{-\nu_{0}(1-\varepsilon)r}$, where $C>0$ is a constant to be chosen later, $\varepsilon>0$ is sufficiently small, $\kappa=\sqrt{\frac{1}{4}g^{2}\rho_{0}^{2}-A_{0}^{2}},\,
\nu=\frac{1}{2}\rho_{0}\sqrt{g^{2}+g'^{2}},\,\nu_{0}=\min\{2\kappa,\nu\}$.~Then for any $\varepsilon>0$, there exist a sufficiently large $r_{\varepsilon}>0$ so that
\begin{eqnarray}
\label{107}
(\tau-\eta)''&=&\frac{1}{4}(g^{2}+g'^{2})\rho^{2}\tau
+\frac{2}{r}(h^{2}B-f^{2}A)-\nu_{0}^{2}(1-\varepsilon)^{2}\eta\notag\\[1mm]
&\geqslant&\nu_{0}^{2}\bigg(1-\frac{\varepsilon}{2}\bigg)^{2}\tau
-\nu_{0}^{2}(1-\varepsilon)^{2}\eta
+\frac{2}{r}(h^{2}B-f^{2}A)\notag\\[1mm]
&\geqslant&\nu_{0}^{2}\bigg(1-\frac{\varepsilon}{2}\bigg)^{2}(\tau-\eta)+I_{1}\notag\\[1mm]
&\geqslant&\nu_{0}^{2}\bigg(1-\frac{\varepsilon}{2}\bigg)^{2}(\tau-\eta),\,r>r_{\varepsilon}.
\end{eqnarray}
Taking $C>0$ be large enough to make
\begin{eqnarray*}
I_{1}&=&\nu_{0}^{2}\left[\left(1-\frac{\varepsilon}{2}\right)^{2}
-\left(1-\varepsilon\right)^{2}\right]\eta
+\frac{2}{r}\left(C_{1}^{2}e^{-2\zeta\left(1-\varepsilon\right)r}A_{0}
-C_{2}^{2}e^{-2\kappa\left(1-\varepsilon\right)r}A_{0}\right)\\[1mm]
&=&\nu_{0}^{2}\left[\left(1-\frac{\varepsilon}{2}\right)^{2}
-\left(1-\varepsilon\right)^{2}\right]Ce^{-\nu_{0}\left(1-\varepsilon\right)r}
+\frac{2}{r}\left(C_{1}^{2}e^{-2\zeta\left(1-\varepsilon\right)r}A_{0}
-C_{2}^{2}e^{-2\kappa\left(1-\varepsilon\right)r}A_{0}\right)>0
\end{eqnarray*}
and $(\tau-\eta)(r_{\varepsilon})=\left\{r_{\varepsilon}
\left[B\left(r_{\varepsilon}\right)-A\left(r_{\varepsilon}\right)\right]
-Ce^{-\nu_{0}\left(1-\varepsilon\right)r_{\varepsilon}}\right\}\leqslant0$, where $h(r)=C_{1}e^{-\zeta(1-\varepsilon)r}$, $C_{1}$ is an arbitrary constant, $f(r)=C_{2}e^{-\kappa(1-\varepsilon)r}$, $C_{2}$ is also an arbitrary constant.~Therefore, according to this and the boundary condition $\left(\tau-\eta\right)(r)\rightarrow0\,\left(r\rightarrow\infty\right)$, we obtain by
applying the maximum principle theorem in $\eqref{107}$ the result $0<\tau\leqslant\eta=Ce^{-\nu_{0}(1-\varepsilon)r},\,r>r_{\varepsilon}$ as expected in $\eqref{05}$.

Then we consider the estimate for $\rho$.~For the new function $T(r)=r\left(\rho-\rho_{0}\right)$, the equation $\eqref{002}$ gives us
\begin{eqnarray*}
T''=\frac{\lambda}{2}\left(\rho+\rho_{0}\right)\rho T
+\frac{1}{2r}\left[f^{2}-\frac{1}{2}r^{2}\left(A-B\right)^{2}\right]\rho.
\end{eqnarray*}

It is seen that $T''=\frac{\lambda}{2}\left(\rho_{0}+\rho_{0}\right)\rho T=\lambda\rho_{0}^{2}T$ as $r\rightarrow\infty$.~Set $\eta(r)=Ce^{-\sqrt{2}\mu_{0}(1-\varepsilon)r}$ with $\varepsilon>0,\,C>0,\,\mu=\sqrt{\frac{\lambda\rho_{0}}{2}},
\,\mu_{0}=\min\{\sqrt{2}\mu,2\kappa,2\nu_{0}\}$.~Hence for any $\varepsilon>0$, there exist a sufficiently large $r_{\varepsilon}>0$ so that
\begin{eqnarray}
\label{109}
\left(T-\eta\right)''&=&\frac{\lambda}{2}\left(\rho+\rho_{0}\right)\rho T
+\frac{1}{2r}\left[f^{2}-\frac{1}{2}r^{2}\left(A-B\right)^{2}\right]\rho
-2\mu^{2}\left(1-\varepsilon\right)^{2}\eta\notag\\[1mm]
&\geqslant&\frac{\lambda}{2}\left(\rho+\rho_{0}\right)\rho\left(T-\eta\right)+I_{2}
\geqslant\frac{\lambda}{2}\left(\rho+\rho_{0}\right)\rho(T-\eta),\,r>r_{\varepsilon},
\end{eqnarray}
where $r_{\varepsilon}>0$ is sufficiently large.~Then we can choose $C>0$ large enough to make
\begin{eqnarray*}
I_{2}&=&\left[2\mu^{2}\left(1-\frac{\varepsilon}{2}\right)^{2}
-2\mu^{2}\left(1-\varepsilon\right)^{2}\right]\eta
+\frac{1}{2r}\left[f^{2}-\frac{1}{2}r^{2}\left(A-B\right)^{2}\right]\rho\\[1mm]
&=&2\mu^{2}\varepsilon\left(1-\frac{3}{4}\varepsilon\right)
Ce^{-\sqrt{2}\mu_{0}\left(1-\varepsilon\right)r}
+\frac{1}{2r}\left[{C_{2}^{2}e^{-2\kappa\left(1-\varepsilon\right)r}}
-\frac{1}{2}{C_{3}^{2}e^{-2\nu_{0}\left(1-\varepsilon\right)r}}\right]\rho_{0}>0
\end{eqnarray*}
and $\left(T-\eta\right)\left(r_{\varepsilon}\right)=\left\{r_{\varepsilon}
\left[\rho(r_{\varepsilon})-\rho_{0}\right]-
Ce^{-\sqrt{2}\mu_{0}\left(1-\varepsilon\right)r_{\varepsilon}}\right\}\leqslant0$ where $f(r)=C_{2}e^{-\kappa(1-\varepsilon)r}$, $C_{2}$ is an arbitrary constant, $B(r)-A(r)=C_{3}r^{-1}e^{-\nu_{0}(1-\varepsilon)r}$, $C_{3}$ is also an arbitrary constant.~Since the boundary condition $(T-\eta)(r)\rightarrow0\,(r\rightarrow\infty)$, using the maximum principle theorem in $\eqref{109}$, then the decay estimate for
$\rho$ near infinity stated in $\eqref{02}$ is established.

Then we study the asymptotic estimates $\eqref{06}$.~Setting
$S(r)=r(\sigma-\sigma_{0})$, for $\eqref{006}$, we obtain that
$S''=2\kappa\sigma_{0}^{2}S$ as $r\rightarrow\infty$.~To get the
estimate for $\sigma(r)$ in $\eqref{06}$, we use the comparison
function $\eta(r)=Ce^{-\sqrt{2}\xi(1-\varepsilon)r}$, where $C>0$ is
a constant to be chosen later, $\varepsilon>0$ is sufficiently
small, $\xi=\sqrt{\kappa}\sigma_{0}$.~Then for any $\varepsilon>0$,
there is a sufficiently large $r_{\varepsilon}>0$ such that
\begin{eqnarray}
\label{111}
\left(S-\eta\right)''&=&\kappa\left(\sigma_{0}+\sigma\right)\sigma S+\frac{2}{r}h^{2}\sigma
-2\xi^{2}\left(1-\varepsilon\right)^{2}\eta\notag\\[1mm]
&\geqslant&\kappa\left(\sigma_{0}+\sigma\right)\sigma\left(S-\eta\right)
+\left[2\xi^{2}\left(1-\frac{\varepsilon}{2}\right)^{2}
-2\xi^{2}\left(1-\varepsilon\right)^{2}\right]\eta\notag\\[1mm]
&\geqslant&\kappa\left(\sigma_{0}+\sigma\right)\sigma\left(S-\eta\right),
\,r>r_{\varepsilon}.
\end{eqnarray}

Let $C>0$ be large enough to make $(S-\eta)(r_{\varepsilon})
=\{r_{\varepsilon}[\sigma(r_{\varepsilon})-\sigma_{0}]
-Ce^{-\sqrt{2}\xi(1-\varepsilon)r_{\varepsilon}}\}\leqslant0$.~Thus, in view of this and the boundary condition $(S-\eta)(r)\rightarrow0\,(r\rightarrow\infty)$, we obtain by applying the maximum principle theorem in $\eqref{111}$ the result $0<S\leqslant\eta=Ce^{-\sqrt{2}\xi(1-\varepsilon)r},\,r>r_{\varepsilon}$ as expected in $\eqref{06}$.

Then the estimate for $A$ in \eqref{03} follows.
The proof of the Theorem \ref{th2} is  completed.
\hfill$\Box$

\end{document}